\documentclass[12pt]{article}
\usepackage[english]{babel}
\usepackage{sfmath}
\usepackage{relsize} 
\usepackage{amssymb}
\usepackage{natbib}
\usepackage[leqno]{mathtools}
\usepackage{upgreek}

\usepackage{amsthm}

\usepackage{xr-hyper}
\usepackage[bookmarks,pdfnewwindow,plainpages=false]{hyperref}

\hypersetup{
colorlinks=true,
linkcolor=black,
citecolor=black,
pdfauthor={Victor Ivrii},
pdftitle={Sharp Spectral Asymptotics for four-dimensional Schr\"odinger operator with a strong  magnetic field. II},
pdfsubject={Sharp Spectral Asymptotics},
pdfkeywords={Microlocal Analysis, Magnetic Schr\"odinger Operator},
}
\newcommand\MW{{\rm {MW}}}

\usepackage{graphicx,color}

\usepackage{eso-pic}
\definecolor{lightgray}{gray}{.90}
\definecolor{verylightgray}{gray}{.93}

\usepackage{subfigure}
\usepackage{placeins}

\usepackage{enumitem}

\newcommand\W{{\rm W}}

\newcommand\Def{\stackrel{\text{def}}{=}}

\newcommand\const{{\rm {const}}}
\newcommand\dist{{\rm {dist}}}

\newcommand\corr{{\rm {corr}}}

\newcommand\Hess{\operatorname{Hess}}

\newcommand\mes{\operatorname{mes}}

\newcommand\rank{\operatorname{rank}}
\renewcommand\Re{\operatorname{Re}}

\newcommand\supp{\operatorname{supp}}
\newcommand\bR{{\mathbb R}}

\newcommand\bZ{{\mathbb Z}}
\newcommand\cA{{\mathcal A}}
\newcommand\cB{{\mathcal B}}

\newcommand\cE{{\mathcal E}}

\newcommand\cQ{{\mathcal Q}}

\newtheorem{theorem}{Theorem}[section]

\newtheorem{corollary}[theorem]{Corollary}
\newtheorem{proposition}[theorem]{Proposition}

\theoremstyle{definition}

\newtheorem{conjecture}{Conjecture}

\theoremstyle{remark}
\newtheorem{remark}[theorem]{Remark}

\numberwithin{equation}{section}

\newenvironment{claim}[1][{\rm(\theequation)}]{\refstepcounter{equation}%
\begin{trivlist}
\item[{\hskip\labelsep#1}]}{\end{trivlist}\addvspace{10pt}}

\newcounter{note}

\newenvironment{claim*}[1]{\medskip
\begin{trivlist}
\item[{\hskip\labelsep#1}]}{\medskip\end{trivlist}}

\newenvironment{phantomequation}[1][]{\refstepcounter{equation}}{}

\renewcommand\subsubsection{\paragraph{\thesubsubsection}\refstepcounter{subsubsection}}

\begin{document}

%\AddToShipoutPicture*{%
%    \AtTextCenter{%
%      \makebox(0,0)[c]{\resizebox{\textwidth}{!}{%
%        \rotatebox{55}{\textsf{\textbf{\color{lightgray} Final}}}}} 
%    }
%  }
%  
%  

\title{\
Sharp spectral asymptotics for generic 4-dimensional Schr\"odinger operator with the strong magnetic field}
\author{Victor Ivrii}
\date{\today}

\maketitle

{\abstract%
I consider 4-dimensional Schr\"odinger operator with the generic non-degenerating magnetic field and for a generic potential I derive spectral asymptotics with the remainder estimate $O(\mu^{-1}h^{-3})$ and the principal part $\asymp h^{-4}$ where $h\ll 1$ is Planck constant and $\mu \gg 1$ is the intensity of the magnetic field. For general potentials remainder estimate $O(\mu^{-1}h^{-3}+\mu^2h^{-2})$ is achieved.
\endabstract}

\section{Introduction}

Sharp spectral asymptotics for multidimensional Magnetic Schr\"odinger were obtained in \cite{IRO4, IRO5} in full- and non-full-rank cases respectively\footnote{\label{foot-1} I mean the rank of magnetic intensity matrix $(F_{jk})$.}. The results, as one could expect from the analysis of 2- and 3-dimensional cases, were rather different.

However there are two problems with these papers: first, the rank of the magnetic intensity matrix there was supposed to be constant which is not necessarily the case even if  only generic magnetic fields are considered; second, while both the maximal rank\footnote{\label{foot-2} I mean that the rank of $(F_{jk})$ is $2\lfloor d/2\rfloor$ at each point.} and microhyperbolicity\footnote{\label{foot-3} Which was assumed in \cite{IRO4}.} conditions  are stable with respect to the small perturbations, they both are not generic in the sense that even in a small but fixed domain   a general Magnetic Schr\"odinger operator is not necessarily approximated by operators satisfying any of these conditions: exactly in the same way as stationary points are not necessarily removable. 

However the analysis of \cite{IRO4} and \cite{IRO6} leads me to the following conjecture for even-dimensional Magnetic Schr\"odinger operators:

\begin{conjecture}\label{conj} As $\mu \le ch^{-1}$ the main part of asymptotics is given by $\cE^\MW$ while 
\begin{enumerate}[label=(\roman*)]

\item For fixed $(g^{jk})$ and general $(V_j), V$ the remainder estimate  is $O(\mu h^{1-d})$;

\item For fixed $(g^{jk})$ and generic $(V_j)$ (i.e. $(V_1,\dots, V_d)\notin {\mathfrak A}_g$ which is nowhere dense closed set) and general $V$ the remainder estimate is  $O(\mu ^{-1} h^{1-d}+\mu^{d/2}h^{-d/2})$ as $(F_{jk})$ has the full rank everywhere and $O(\mu ^{-1/2} h^{1-d}+\mu^{d/2}h^{-d/2})$ otherwise;

\item For fixed $(g^{jk})$ and generic $(V_j,V)$ (i.e. $(V_1,\dots, V_d;V)\notin {\mathfrak B}_g$ which is nowhere dense closed set)  the remainder estimate is  $O(\mu ^{-1} h^{1-d})$ as $(F_{jk})$ has the full rank everywhere and $O(\mu ^{-1/2} h^{1-d})$ otherwise.
\end{enumerate}
\end{conjecture}

While (i) is trivial, (ii) and (iii) are rather difficult and
my goal is rather limited: to prove them  as $d=4$. For $d=2$ it was done in \cite{IRO6, IRO16}). Further, in the small vicinity of the set 
$\{x: \rank (F_{jk})(x)\le 2\}$ it was done in \cite{IRO8}. So, I will need to investigate the case of $\rank (F_{jk})=1$ at every point.

So, operator in question is
\begin{equation}
A= {\frac 1 2}\Bigl(\sum_{j,k}P_jg^{jk}(x)P_k -V\Bigr),\qquad P_j=D_j-\mu V_j
\label{0-1}
\end{equation}
with smooth\footnote{\label{foot-4} Smooth means either infinitely smooth or belonging to $C^K$ with large enough $K$.}  symmetric positive definite matrix $(g^{jk})$ and smooth real-valued potentials  $(V_1,\dots, V_d;V)$ and $\mu \gg1$, $h\ll 1$. Assuming that $A$ is self-adjoint, let $E(\tau)$ be the spectral projector of $A$ and $e(x,y,\tau)$ be its Schwartz' kernel. These assumptions are fulfilled during the whole article.

Magnetic field is characterized by a skew-symmetric matrix $(F_{jk})$, $F_{jk}=\partial_jV_k-\partial_kV_j$ and more precisely  by $(F^j_k)=(g^{j*})(F_{*k})$ and its eigenvalues $\pm if_j$, $f_j\ge 0$. As $d=4$ these are $f_1$ and $f_2$. 

It is proven \cite{Ma} that 

\begin{claim}\label{0-2} For generic $(V_1,\dots,V_4)$ $f_1$ and $f_2$ do not vanish simultaneously.
\end{claim}

\begin{remark}\label{rem-0-1}  (i) If  one of $f_1,f_2$ vanishes then (for generic $(V_1,\dots,V_4)$ locally situation of \cite{IRO8} occurs and the remainder estimate is $O(\mu^{-1/2}h^{-3})$ for generic $V$  and 
$O(\mu^{-1/2}h^{-3}+\mu^2h^{-2})$ for general $V$ (detailed assumptions see in 
 \cite{IRO8}) and I exclude this case from the analysis, assuming that
 \begin{equation}
f_1 \ge \epsilon_0,\ f_2\ge \epsilon_0.
\label{0-3}
\end{equation}

\medskip
\noindent
(ii) As condition (\ref{0-4}) is fulfilled and $\mu h \ge c$ then the main part of asymptotics is $0$ and the remainder estimate is $O(\mu^{-s})$ with an arbitrarily large exponent $s$ (see f.e. \cite{IRO3}) and therefore in what follows I assume that
\begin{equation}
1\le \mu \le ch^{-1}.
\label{0-4}
\end{equation}
\end{remark}

My goal is to prove two following theorems:

\begin{theorem}\label{thm-0-2}
Let $(g^{jk})$ be fixed and $(V_1,\dots,V_4)$ be generic. Further, let conditions $(\ref{0-3}),(\ref{0-4})$ be fulfilled and $\psi$ be a smooth function. Then
\begin{equation}
|\int \bigl( e (x,x,0)-\cE^\MW (x,0)\bigr)-\cE_\corr^\MW (x,0)\bigr)\psi (x)\,dx| \le C\mu^{-1}h^{-3}+ C\mu^2 h^{-2}
\label{0-5}
\end{equation}
where
\begin{multline}
\cE^\MW (x,\tau)= \\(2\pi)^{-2}\mu^2h^{-2}\sum _{(m,n)\in \bZ^{+\,2} } \theta \bigl(2\tau+V- (2m+1)\mu h f_1 - (2n+1)\mu hf_2\bigr) f_1 f_2 \sqrt g
\label{0-6}
\end{multline}
is Magnetic Weyl Expression, $g=\det (g^{jk})^{-1}$, $\cE_\corr^\MW$ is a correction term, defined by $(\ref{4-23})$.
\end{theorem}

\begin{theorem}\label{thm-0-3}
Let $(g^{jk})$ be fixed and $(V_1,\dots,V_4;V)$ be generic. Further, let conditions $(\ref{0-3}),(\ref{0-4})$ be fulfilled and $\psi$ be a smooth function. Then
\begin{equation}
|\int \bigl( e (x,x,0)-\cE^\MW (x,0)\bigr)\psi (x)\,dx| \le C\mu^{-1}h^{-3}.
\label{0-7}
\end{equation}
\end{theorem}

\begin{remark}\label{rem-0-4} (i) More precise conditions describing conditions for  $(V_1,\dots, V_d)$ in Theorem \ref{thm-0-2}  and for $(V_1,\dots, V_d;V)$ in Theorem \ref{thm-0-3} will be formulated below in section \ref{sect-4}.

\medskip
\noindent
(ii) I was able to prove that one can skip $\cE_\corr^\MW (x,0)$ and preserve remainder estimate (\ref{0-5}) unless $h^{-1/3+\delta}\le \mu \le h^{-1/3-\delta}$; in the latter case the remainder estimate (\ref{0-6}) should be replaced by $h^{-8/3-\delta}$; however I suspect that one can always skip $\cE_\corr^\MW (x,0)$ and preserve remainder estimate (\ref{0-5}).

\medskip
\noindent
(iii)  With the exception of section \ref{sect-5} I assume that 
 \begin{equation}
V\ge \epsilon_0.
\label{0-8}
\end{equation}
In section \ref{sect-5} I will get rid off this condition.
\end{remark}

\paragraph{Plan of the paper.} Section \ref{sect-1} is devoted to the geometry (discussion of what is the generic case) and the preliminary analysis in the cases when results of \cite{IRO4} imply  theorem \ref{thm-0-3} immediately. 

In section \ref{sect-2} I tackle the weak magnetic field case ($\mu\le h^{-\delta_0}$ with sufficiently small exponent $\delta_0>0$)  proving the standard Weyl formula with the remainder $O(\mu^{-1}h^{-3})$.

Sections   \ref{sect-3}, \ref{sect-4} are devoted to the case when $ h^{-\delta_0}\le \mu \le ch^{-1}$ with arbitrarily small exponent ${\bar\delta}>0$ and one can reduce operator to the microlocal canonical form. More precisely, in section \ref{sect-3} I prove asymptotics with the announced remainder estimates but with the implicitly given main part
\begin{equation}
h^{-1}\sum_{\iota}\int_{-\infty}^0 \Bigl( F_{t\to h^{-1}\tau} {\bar\chi}_{T_\iota}(t) \Gamma uQ_{\iota y}^t\Bigr)\,d\tau
\label{0-9}
\end{equation}
where $u$ is the Schwartz' kernel of the propagator $e^{ih^{-1}tA}$ and $Q_\iota$ are pseudo-differential partition elements (see details in my previous papers); in this proof estimate of
\begin{equation}
|\Bigl( F_{t\to h^{-1}\tau} {\bar\chi}_{T_\iota}(t) \Gamma uQ_{\iota y}^t\Bigr)|
\label{0-10}
\end{equation}
plays the crucial role.

In section \ref{sect-4} I replace (\ref{0-10}) by the standard Magnetic Weyl formula and estimate an error arising from this.

Finally, in section  \ref{sect-5} I consider the case when condition (\ref{0-8}) is violated.

\section{Geometry and Preliminary Analysis}\label{sect-1}

\subsection{Geometry}\label{sect-1-1}

\begin{proposition}\label{prop-1-1} Let either $g^{jk}$ be fixed and then $V_j$ be generic, or,  alternatively, $V_j$ be fixed and then $g^{jk}$ be generic. Further, let $(\ref{0-3})$ be fulfilled. Then
\begin{claim}\label{1-1}
$\Sigma\Def\{x:f_1=f_2\}$  is a smooth 2-dimensional manifold: $\Sigma=\{x:v_1=v_2=0\}$ with appropriate smooth functions $v_1$ and $v_2$ such that $\nabla v_1$ and $\nabla v_2$ linearly independent at any point of $\Sigma$ and 
\begin{equation}
\dist (x,\Sigma)\asymp|f_1-f_2| = 2  (v_1^2+v_2^2)^{1/2};
\label{1-2}
\end{equation}
furthermore, $f_1+f_2$ and $(f_2-f_1)^2$ are  smooth functions;
\end{claim} 

\begin{claim}\label{1-3} Consider symplectic form corresponding to $F_{jk}$
\begin{equation}\label{1-4}
\omega={\frac 1 2}\sum_{jk}F_{jk}dx_j\wedge dx_k;
\end{equation}
then  $\omega|_\Sigma$ is either non-degenerate or it is generic degenerate\footnote{\label{foot-5}  Exactly as $\Sigma$ was  in  \cite{IRO6}.} i.e. $\Theta_1\Def\{x\in \Sigma,\ \{v_1,v_2\}=0\}$ is a submanifold of dimension $1$ and  $\nabla_\Sigma \{v_1,v_2\}$ is disjoint from 0 at $\Theta_1$;
\end{claim}

\begin{claim}\label{1-5} Function $f_1f_2^{-1}$ has only non-degenerate critical points outside of $\Sigma$ and near $\Sigma$ $|\nabla (f_1f_2^{-1})|$ is disjoint from $0$.
\end{claim}
\end{proposition}

\begin{proof} I leave an easy proof to the reader.
\end{proof}

\subsection{Microhyperbolicity Condition}
\label{sect-1-2}

\subsubsection{}
\label{sect-1-2-1}
First, let me discuss the case when magnetic intensities $f_1$ and $f_2$ are disjoint. Then they are smooth. It still does not exclude third-order resonances 
$f_1=2f_2$ and $f_2=2f_1$\,\footnote{\label{foot-6} Higher-order resonances are not a problem according to \cite{IRO4}, at least under microhyperbolicity condition.}.

Then there are two cyclotron movements and the drift with the velocity $O(\mu^{-1})$ which I want to calculate. Let 
\begin{equation}
p_j=\xi_j - \mu V_j.
\label{1-6}
\end{equation}
Note that 
\begin{align}
&x'_j \Def x_j-\mu^{-1}\sum_l\phi^{jl}p_l,\label{1-7}\\
&(\phi^{**})\Def (F_{**})^{-1}\label{1-8}
\end{align}
satisfy 
\begin{equation}
\{p_k,x_j'\}=-\mu^{-1}\sum_l \{p_k,\phi^{jl}\}p_l
\label{1-9}
\end{equation}
and therefore 
\begin{multline}
\Bigl\{{\frac 1 2}\sum_{k,m} g^{km}p_kp_m, x'_j\Bigr\}= {\frac 1 2}\sum_{k,m} \{g^{km},x'_j\} p_kp_m -\mu^{-1}\sum_{k,m,l}  g^{km} \{p_k,\phi^{jl}\}p_mp_l=\\
\mu^{-1} \sum_{k,m,l}\Bigl(   {\frac 1 2}\phi^{jl}\{p_l,g^{km}\} -\{p_l,\phi^{jk}\}g^{ml}
\Bigr)p_kp_m
\label{1-10}
\end{multline}

\begin{proposition}
Let $f_1\ne f_2$ in $\Omega$; then one can correct 
\begin{equation}
x_j'\mapsto x_j''\Def x_j' -{\frac 1 2}\mu^{-1}\sum_{k,m}\beta^{jkm}p_kp_m
\label{1-11}
\end{equation}
so that 
\begin{equation}
\Bigl\{{\frac 1 2}\sum_{k,m} g^{km}p_kp_m, x''_j\Bigr\}= 
\mu^{-1}\Bigl(\{f_1,x_j\}_F b_1 + \{f_2,x_j\}_F b_2\Bigr)+ O(\mu^{-2})
\label{1-12}
\end{equation}
where in the Birkhoff normal form $a=f_1 b_1+f_2b_2+\dots$ and 
\begin{equation}
\{w_1,w_2\}_F\Def \sum_{k,l}\phi^{kl}\partial_k w_1\cdot \partial_l w_2 
\label{1-13}
\end{equation}
are Poisson brackets associated with symplectic form $\omega=\omega_F$.
\end{proposition}

\subsubsection{}\label{sect-1-2-2}
\emph{Microhyperbolicity condition\/} of \cite{IRO4} then means for $d=4$ and $f_1\ne f_2$  that
\begin{multline}
|V-(2p+1)f_1\mu h - (2n+1)f_2\mu h|\le \epsilon \implies\\
|\langle \ell,\nabla \bigl(V-(2p+1)f_1\mu h - (2n+1)f_2\mu h\bigr)\rangle|\ge \epsilon
\label{1-14}
\end{multline}
where $\ell$ ($|\ell|\asymp 1$) is the \emph{microhyperbolicity direction\/} (at the given point $x$) and $(n,p)$ are \emph{magnetic indices\/}. However there actually were two microhyperbolicity conditions in \cite{IRO4}: in the weaker condition $\ell$ could depend not only on $x$ but also on magnetic indices while in the stronger one it was assumed that $\ell$ is the same for all pairs $(n,p)$ connected by the third order resonance. Respectively, the remainder estimates derived in \cite{IRO4} were $O(h^{1-d})$ under weaker condition and $O(\mu^{-1}h^{1-d})$ under stronger one.

In the relatively simple case of $d=4$ the resonance means $kf_1=lf_2$ where $k,l\in \bZ^+$ are coprimes, $(k+l)$ is the \emph{order  of the resonance\/}. So under the stronger microhyperbolicity condition if either $2f_1-f_2=0$ or $f_1-2f_2=0$ at $x$ then $\ell$ should not depend on $(n,p)$ satisfying the the left-hand inequality in (\ref{1-14}).

\begin{remark}\label{rem-1-3} Case $d=4$ is relatively simple  because there are only two magnetic intensities $f_1$ and $f_2$ and thus it is impossible to have 
collisions between two second-order resonances (of the type $f_j=f_k$), two third order resonances (of the types $f_j=2f_k$ or $f_j=f_k+f_l$) or the second and the third order resonances.

Further, as $|kf_1-lf_2|\le \epsilon$ there could be other resonances only  of order $m(\epsilon)$ or higher; $m(\epsilon)\to +\infty$ as $\epsilon\to +0$.
\end{remark}

So, as $d=4$ and $f_1$ is disjoint from $f_2$ and $\mu h\le \epsilon_1$ (where $\epsilon_1>0$ depends on $\epsilon$ in the microhyperbolicity condition) both microhyperbolicity conditions are equivalent  to
\begin{equation}
| \nabla\bigl( \alpha\log f_1+(1-\alpha)\log f_2- \log V\bigr)|\ge \epsilon\qquad \forall \alpha: 0\le \alpha\le 1.
\label{1-15}
%\label{1-22}
\end{equation}
On the other hand, as $\mu h\ge \epsilon_1$ condition (\ref{1-15}) should be checked only at points where \emph{ellipticity condition\/}
\begin{equation}
|V-(2p+1)f_1\mu h - (2n+1)f_2\mu h|\ge \epsilon  \qquad \forall (p,n)\in \bZ^{+\,2}
\label{1-16}
\end{equation}
is violated.

\begin{remark}\label{rem-1-4} I remind that under condition (\ref{1-16}) on $\supp\psi$ the remainder estimate is $O(\mu^{-s}h^s)$ with an arbitrarily large exponent $s$ and the similar result would hold in any dimension provided $\rank (F_{jk})=d$ at each point.
\end{remark}

So, main theorems \ref{IRO4-thm-0-3}, \ref{IRO4-thm-0-5}, \ref{IRO4-thm-0-7}   of \cite{IRO4} imply 

\begin{proposition}\label{prop-1-5} Let $d=4$ and  conditions $(\ref{0-8})$, $f_1\ne f_2$ and  $(\ref{1-15})$ be fulfilled at $\supp\psi$.  Then  the remainder is $O(\mu^{-1}h^{-3})$ while the main part of asymptotics is given by  the Magnetic Weyl formula.
\end{proposition}

\subsubsection{}\label{sect-1-2-3} This  leaves us with bad points where condition (\ref{1-15}) is violated.
So let us consider set $\Lambda_\alpha$ of critical points of 
\begin{equation}
\phi_\alpha\Def  \alpha\log f_1+(1-\alpha)\log f_2- \log V.
\label{1-17}
\end{equation}

\begin{proposition}\label{prop-1-6}  Let $g^{jk}$, $V_j$ be fixed and let 
$V$ be generic. Then 

\medskip
\noindent
{\rm (i)} $\Lambda_\alpha$ is a finite set, continuously depending on $\alpha\in [0,1]$;

\medskip
\noindent
{\rm (ii)} There exist $0<\alpha_1 < \alpha_2< \dots <\alpha_J<1$ such that for $\alpha \notin \{\alpha_1,\dots, \alpha_J\}$  function $\phi_\alpha$ has non-degenerate critical points, while for $\alpha=\alpha_j$ function  
$\phi_\alpha$ has also isolated critical points but in one of them $\rank \Hess \phi_\alpha=3$;

\medskip
\noindent
{\rm (iii)} As ${\bar\alpha}\in \{\alpha_1,\dots, \alpha_J\}$ and ${\bar x}$ is a critical point of $\phi_{\bar\alpha}$ with $\rank \Hess \phi_{\bar\alpha}({\bar x})=3$, in an appropriate (smooth) coordinate system in vicinity of ${\bar x}=0$
\begin{equation}
\phi_\alpha =   {\frac 1 3}x_1^3 -(\alpha -{\bar\alpha})x_1 
\pm_2 {\frac 1 2}x_2^2 \pm_3 {\frac 1 2}x_3^2  \pm_4 {\frac 1 2}x_4^2
\label{1-18}
\end{equation}
with independent signs $\pm_k$.
\end{proposition}

\begin{proof} An easy proof is left to the reader.
\end{proof}

\begin{remark}\label{rem-1-7} (i)  Assume that at some point ${\bar x}$
\begin{equation}
 |\nabla( f_1f_2^{-1}) |\ge \epsilon_0.
\label{1-19}
\end{equation}
Then in the vicinity of ${\bar x}$
\begin{equation}
|\nabla \phi_\alpha |\ge \epsilon |\alpha -{\bar\alpha}(x)|\qquad \forall \alpha: 0\le \alpha\le 1
\label{1-20}
\end{equation}
with smooth function ${\bar\alpha}(x)$.

\medskip
\noindent
(ii) On the other hand, if ${\bar x}$ is an isolated critical point of $f_1f_2^{-1}$ then as $V$ is generic  $|\nabla \phi_\alpha|\ge \epsilon$ for all $\alpha \in [0,1]$ in the vicinity of ${\bar x}$.
\end{remark}

\subsection{Analysis near $\Sigma$: Geometry}
\label{sect-1-3}

\subsubsection{}
\label{sect-1-3-1}

Consider now $\Sigma=\{x: f_1 =f_2\}$ assuming that (\ref{1-1}) holds; then the microhyperbolicity condition at point $x\in \Sigma$ means exactly that
\begin{equation}
{\frac 1 2}\ell \bigl((f_1+f_2)V^{-1}\bigr)\ge 
\bigl( (\ell v _1)^2+(\ell v_2)^2\bigr)^{1/2}V^{-1} +\epsilon
\label{1-21}
\end{equation}
where $\ell$ again is the microhyperbolicity direction, or, equivalently, 
\begin{equation}
|{\frac 1 2}\nabla  \bigl((f_1+f_2)V^{-1}\bigr) -(\beta_1 \nabla v_1 +\beta_2\nabla v_2)V^{-1}|\ge 
\epsilon \qquad \forall \beta =(\beta_1,\beta_2)\in \bR^2: |\beta|\le 1.
\label{1-22}
\end{equation}

\begin{remark}\label{rem-1-8} Note that the microhyperbolicity condition could be violated only at stationary points of $f_1V^{-1}|_\Sigma$ and only at those of them where 
\begin{equation}
{\frac 1 2}\nabla  \bigl((f_1+f_2)V^{-1}\bigr) =(\beta_1 \nabla v_1 +\beta_2\nabla v_2)V^{-1}
\label{1-23}
\end{equation}
with $\beta_1^2+\beta_2^2\le 1$.
\end{remark}

Let $\Sigma_0$ be the set of stationary points of $f_1V^{-1}|_\Sigma$; at each point of $\Sigma_0$ decomposition (\ref{1-23}) holds; let us denote by $\Sigma_0^+$, $\Sigma_0^-$, $\Sigma_0^0$ the subsets of $\Sigma_0$ where $\beta_1^2+\beta_2^2>1$, $\beta_1^2+\beta_2^2<1$ and $\beta_1^2+\beta_2^2=1$ respectively. Then the microhyperbolicity condition holds at $\Sigma_0^+$.

So, main theorems \ref{IRO4-thm-0-3}, \ref{IRO4-thm-0-5}, \ref{IRO4-thm-0-7}   of \cite{IRO4} imply 

\begin{proposition}\label{prop-1-9} Let $d=4$ and  conditions $(\ref{1-1})$ and $(\ref{0-8})$  be fulfilled.   Let $\psi$ be supported in the small vicinity of ${\bar x}\in \Sigma_0^+$. Then  the remainder is $O(\mu^{-1}h^{-3})$ while the main part of asymptotics is given by  magnetic Weyl formula.
\end{proposition}

\subsubsection{}
\label{sect-1-3-2}

One can prove easily two following propositions:

\begin{proposition}\label{prop-1-10} Let $d=4$ and  condition $(\ref{1-1})$  be fulfilled.  Let $V$ be generic satisfying $(\ref{0-8})$. Then
\begin{claim}\label{1-24} $\Sigma_0$ consists of the finite number of non-degenerate stationary points of $f_1V^{-1}|_\Sigma$; these points are generic;\end{claim}
\begin{claim}\label{1-25}
 $0<\beta_1^2+\beta_2^2<1$ at each point of $\Sigma_0 \setminus \Sigma_0^+$.
\end{claim}
\end{proposition}

\begin{proposition}\label{prop-1-11} Let $d=4$ and  condition $(\ref{1-1})$  be fulfilled. Let ${\bar x}\in \Sigma_0^-$ with $0<\beta_1^2+\beta_2^2<1$ in decomposition $(\ref{1-23})$. Then in the vicinity of ${\bar x}$ 
\begin{equation}
\Lambda \Def  \bigl\{x: f_1\ne f_2 \text{ and } \exists \alpha\in [0,1]:\quad \nabla\phi_\alpha=0\bigr\}\cup \{{\bar x}\}
\label{1-26}
\end{equation}
is a smooth curve  passing through ${\bar x}$ and transversal to $\Sigma$; moreover  
\begin{equation}
{\bar\alpha}({\bar x}) = {\frac 1 2}\bigl(1\pm (\beta_1^2+\beta_2^2)^{1/2}\bigr)
\label{1-27}
\end{equation}
\end{proposition}

\section{Weak Magnetic Field Case}
\label{sect-2}

In this section I am going to prove the remainder estimate $O(\mu^{-1}h^{-3})$ for general $V$ as
\begin{equation}
c\le \mu \le h^{-\delta_0}
\label{2-1}
\end{equation}
with small exponent $\delta_0>0$.

\setcounter{subsection}{1}
\setcounter{subsubsection}{0} 
\subsubsection{}\label{sect-2-1-1} Assume first that condition (\ref{1-19}) is fulfilled.

\begin{proposition}\label{prop-2-1} Let conditions $f_1\ne f_2$ and  $(\ref{1-19})$ be fulfilled at $\supp\psi$.  Then under condition $(\ref{2-1})$ the remainder is $O(\mu^{-1}h^{-3})$ while the main part of asymptotics is given by  Weyl formula.
\end{proposition}

\begin{proof}[Proof, Part I] (i) Assume first that there are no  resonances of order lesser than $M$ on $\supp \psi$:
\begin{equation}
|kf_1-lf_2|\ge \epsilon,\qquad \forall k,l\in \bZ^+: k+l\le M.
\label{2-2}
\end{equation}
Then one can reduce operator to the normal form without cubic and ``unbalanced'' fourth order terms and then along Hamiltonian trajectories 
\begin{equation}
{\frac {d\ }{dt}}x_j= {\frac {\partial\ }{\partial\xi_j}}a,\qquad
{\frac {d\ }{dt}}\xi_j= -{\frac {\partial\ }{\partial x_j}}a
\label{2-3}
\end{equation}
of
\begin{equation}
a(x,\xi)\Def {\frac 1 2}\Bigl(\sum_{j,k}p_jg^{jk}(x)p_k -V\Bigr),\qquad p_j=\xi_j-\mu V_j
\label{2-4}
\end{equation}
the following relations hold:
\begin{equation}
{\frac {d\ }{dt}}b_j=\{a , b_j\}= O(\mu^{-2})
\label{2-5}
\end{equation}
where $b_j = {\frac 1 2} (f_j \circ \Psi) (x_j^2+\xi_j^2)+O(\mu^{-2})$ are quadratic terms in the Birkhoff normal form; see details in \cite{IRO4} or below.

Therefore in the ``corrected'' coordinates $x''$  defined by (\ref{1-11}) along Hamiltonian trajectories 
\begin{equation}
{\frac {d\ }{dt}}x''_j= \mu^{-1}\bigl(\beta_{j1}(x'') b_1(x,\xi) + \beta_{j2}(x'')b_2(x,\xi)\bigr)+O(\mu^{-3})
\label{2-6}
\end{equation}
and 
\begin{equation}
|{\frac {d^2\ }{dt^2}}x''|\le C \mu^{-1}|{\frac {d\ }{dt}}x''|+C(\mu^{-4}).
\label{2-7}
\end{equation}
Thus as  
\begin{equation}
\mu |{\frac {d\ }{dt}}x''|\asymp \rho\ge C\mu^{-2}
\label{2-8}
%\label{2-6}
\end{equation}
then for time $T_1=\epsilon \mu $  this relation (\ref{2-8}) is retained and    also variation of $dx''/dt$ would be less than $C\epsilon \rho \mu^{-1}$. 

Then there is a fixed direction $\ell$, $|\ell |\asymp 1$ such that for $|t|\le T_1$
\begin{equation}
\langle\ell,{\frac {d\ }{dt}}x''\rangle \ge \epsilon_0 \mu^{-1}\rho;
\label{2-9}
%\label{2-8}
\end{equation} 
$\ell$ is a sort of the microhyperbolicity direction. Without any loss of the generality one can assume that $\ell=(1,0,0,0)$ and therefore  the shift  of $x''_1$ for time $T$ is exactly of magnitude $\rho \mu^{-1}T$.

According to the logarithmic uncertainty principle this shift is microlocally observable as 
\begin{equation}
\rho \mu^{-1}T \times \rho \ge Ch|\log h|
\label{2-10}
\end{equation}
Plugging $T=T_0=\epsilon\mu^{-1}$ one gets $\rho^2\ge C\mu^2h|\log h|$  which would hold for $\rho\ge C\mu^{-2}$ as $\mu\le h^{-\delta}$. 

Therefore in this case the contribution of zone  (\ref{2-6}) to the remainder does not exceed 
\begin{equation*}
Ch^{-3} \rho \times T_1^{-1} \asymp C\mu^{-1}h^{-3}\rho
\end{equation*}
where factor $\rho h^{-3}$ is due to the estimate of 
\begin{equation}
|F_{t\to h^{-1}\tau}\bigl({\bar\chi}_{T_0}(t)\Gamma Q u\bigr)|\le C\rho h^{-3},
\label{2-11}
\end{equation}
$Q$ is the cut-off operator in zone (\ref{2-8}), $u$ is the Schwartz kernel of $e^{-ih^{-1}tA}$ and other  notations of \cite{IRO3,IRO8} are used.

To prove estimate (\ref{2-11}) let us make $\rho$-admissible partition in $\xi$; then one can prove easily that the contribution of each such element to  the left hand expression does not exceed $C\rho^3h^{-3}$ because the propagation speed with respect to $x$ is $\asymp 1$ as $|t|\le T_0$ and therefore one can trade $T_0$ to ${\bar T}=Ch\rho^{-1}|\log h|$ in the left hand expression of (\ref{2-11}) with a negligible error and in the estimate even to ${\bar T}=C\rho^{-1}h$ (one can prove easily by the rescaling). On the other hand, the number of the partition elements for each $x$ (and therefore in its vicinity) so that $hD_t-A$ is not elliptic and (\ref{2-8}) holds is obviously $O(\rho^{-2})$; so the left-hand expression of (\ref{2-11}) does not exceed $C\rho^3h^{-3}\times \rho^{-2}$.

Summation of $O(\mu^{-1}h^{-3})$ with respect to $\rho$ results in $C\mu^{-1}h^{-3}$ and therefore I have proven that

\begin{claim}\label{2-12} In frames of proposition and (\ref{2-2}) the total contribution of zone 
\begin{equation}
\bigl\{ \mu |{\frac {d\ }{dt}}x''|\ge {\bar\rho}\Def  C\mu^{-2}\bigr\}
\label{2-13}
\end{equation}
to the remainder is $O(\mu^{-1}h^{-3})$.\end{claim}

On the other hand, contribution of zone 
\begin{equation}
\bigl\{ \mu |{\frac {d\ }{dt}}x''|\le {\bar\rho}=C\mu^{-2}\bigr\}
\label{2-14}
\end{equation}
to the remainder does not exceed 
$C{\bar\rho}h^{-3}\times T_0^{-1} =C\mu^{-1}h^{-3}$ where for cut-off operator in this zone estimate (\ref{2-11}) holds with $\rho={\bar\rho}$.

So, the remainder does not exceed $C\mu^{-1}h^{-3}$ while the main part of the asymptotics is given by the standard Tauberian formula (\ref{0-9}) with $T=Ch|\log h|$. One can easily rewrite (\ref{0-9}) with  $T=Ch|\log h|$ as
\begin{equation*}
\int \cE^\W (x,0)\psi(x)\,dx + O(\mu^2h^{-2}).
\end{equation*}

Furthermore, ne can replace here $\cE^\W (x,0)$ by $\cE^\MW (x,0)$ with the same error. The proof is  standard, easy and left to the reader.\end{proof}

\subsubsection{}\label{sect-2-1-2} Let us allow resonances.

\begin{proof}[Proof, Part II]\label{pf-2-1-II}
(i) Let us consider now the case when $(kf_1-lf_2)$ is not disjoint from 0. I will analyze the third order resonance which is worst case scenario leaving the easier case of higher order resonances to the reader.

So, now $(f_1-2f_2)$ is not disjoint from $0$ (but then one can assume that $(kf_1-lf_2)$ is as $l\ne 2k$, $k\le M$). Then one can  reduce a operator to a (pre)canonical form 
\begin{equation}
{\frac 1 2}\Bigl(f_1(x')Z_1^*Z_1 + f_2(x')Z_2^*Z_2 + \mu^{-1}\Re \beta (x') Z_1^*Z_2^2\Bigr) + O(\mu^{-2})
\label{2-15}
%\label{2-32}
\end{equation}
with
\begin{align}
&\{Z_j,Z_k\}\equiv 0, \qquad  \{Z_j^*,Z_k\}\equiv 2\mu \delta_{jk}  &&\mod O(\mu^{-s}).\label{2-16}\\
&\{Z_j,x'_k\}\equiv \{Z_j,^*x'_k\}\equiv 0\qquad &&\mod O(\mu^{-s}),\label{2-17}\\
&\{x'_j,x'_k\}\equiv \mu^{-1}\phi^{jk}(x')\qquad &&\mod O(\mu^{-2})\label{2-18}\\
\intertext{where}
&x'\equiv x\qquad \mod O(\mu^{-1}).\label{2-19}
\end{align}
Due to (\ref{1-19}) one can assume without any loss of the generality that 
\begin{equation}
 f_1-2f_2=x_1.
%\label{2-36}
\label{2-20}
\end{equation}

\medskip
\noindent 
(ii) Now let us consider elements of the partition with 
\begin{equation}
 |x_1| \asymp \gamma, \qquad    \gamma \ge C\mu^{-1}.
%\label{2-36}
\label{2-21}
\end{equation}
Then one can get rid off the cubic term in (\ref{2-15}) by means of the transformation with the generating function 
$\gamma^{-1}\mu^{-2}\Re \beta Z_1^*Z_2^2$
with $\gamma$-admissible function $\beta$ (all other cubic terms are ``regular'' and one can get rid off them in the regular way), leading to the error which, as one can easily check, is the sum of terms of the types $\gamma^{-1}\mu^{-2}\Re \beta Z_j^*Z_j Z_2^*Z_2$ and also some smaller terms; one can continue this process getting rid of all ``unbalanced'' terms up to  order $M$. Also one can introduce corrected $x''$ so that
\begin{equation}
{\frac {d\ }{dt}}x''_j =\mu^{-1} \bigl(\beta_{j1}(x'') b_1 + \beta_{j2}(x'') b_2 \bigr) + O(\mu^{-3}\gamma^{-1})
\label{2-22}
\end{equation}
(compare with (\ref{2-6}); however this would be slightly short of what is needed and one can improve a term $O(\mu^{-3}\gamma^{-1})$ in (\ref{2-22}) to \begin{equation*}
\sum_{k+l=2}\beta_{kl}(x'')b_1^kb_2^l+ O(\mu^{-3-\kappa}\gamma^{-1-\kappa})
\end{equation*}
with $\kappa>0$. Then the contribution of the zone
\begin{equation*}
\bigl\{|x_1|\asymp \gamma, \ \mu |{\frac {d\ }{dt}}x''|\le C\mu^{-2-\kappa}\gamma^{-1-\kappa}\bigr\}
\end{equation*}
to the remainder does not exceed $C\mu h^{-3}\gamma \times \mu^{-2-\kappa}\gamma^{-1-\kappa}\asymp C\mu^{-1-\kappa}h^{-3}\gamma^{-\kappa}$ where the first factor $\gamma$ is the measure. Proof is similar to one in Part I, with estimate (\ref{2-11}) replaced by
\begin{equation}
|F_{t\to h^{-1}\tau}\bigl({\bar\chi}_{T_0}(t)\Gamma Q u\bigr)|\le C\rho\gamma  h^{-3},
\label{2-23}
\end{equation}
now $Q$ is the cut-off operator in zone (\ref{2-8}) intersected with $\{|x_1|\asymp \gamma\}$. Obviously  summation with respect to $\gamma \ge C\mu^{-1}$ results in $O(\mu^{-1}h^{-3})$.

So one needs to consider the zone where
\begin{equation}
\mu |{\frac {d\ }{dt}}x'' |\asymp \rho \ge C\mu^{-2-\kappa}\gamma^{-\kappa},
\label{2-24}
\end{equation}
which implies that 
\begin{equation}
|\mu {\frac {d\ }{dt}}\bigl(\beta_{j1}(x'') b_1 + \beta_{j2}(x'') b_2 \bigr) | \le C\rho + C\mu^{-K}\gamma^{-K}
\label{2-25}
\end{equation}
and therefore one can take 
\begin{equation}
T_1= \epsilon  \min (\mu\gamma \rho^{-1},\mu, \rho \mu^K\gamma^K)
\label{2-26}
\end{equation}
where the restriction $T\le \epsilon \gamma \rho^{-1}$ preserves the magnitude of $|x_1|$.

Therefore the contribution of the zone 
\begin{equation*}
\bigl\{|x_1|\asymp \gamma, \ \mu |{\frac {d\ }{dt}}x''|\asymp \rho \bigr\}
\end{equation*}
with $\gamma\ge C\mu^{-1}$, $\rho \ge C\mu^{-2-\kappa}\gamma^{-\kappa}$ to the remainder does not exceed 
\begin{equation}
C h^{-3}\rho \gamma \Bigl(\mu^{-1}\gamma ^{-1}\rho +1 + \mu^{-K} \gamma^{-K}\Bigr)\asymp C\mu^{-1}h^{-3}\Bigl(\rho^2 +\rho\gamma + \mu^{1-K}\gamma^{1-K}\Bigr).
\label{2-27}
\end{equation}
Obviously summation with respect to $\rho,\gamma$ results in $C\mu^{-1}h^{-3}|\log \mu|$.

\medskip
\noindent
(iii) Now I want to improve this estimate getting rid off $\log \mu$ factor. Note first that  the second term in (\ref{2-25}) sums to $O(\mu^{-1}h^{-3})$, while the first and third terms sum to $O(\mu^{-1}h^{-3})$ over zones complemental to $\{\rho \ge \max(\gamma,|\log \mu|^{-1/2})\}$ and $\{\gamma\le \min (\rho, \mu^{-1+\kappa})\}$ respectively.

However in the first zone  one can take $T_1=\epsilon \mu \gamma^{1-\kappa}$ since in an appropriate time direction $|x_1|\gtrsim \gamma$ for this time. Then an extra factor $\gamma^\kappa$ in the remainder estimate prevents appearance of the logarithmic factor.

Let us introduce $W(x')=Vf_1^{-1}|_{x_1=0}$ and a scaling function $\zeta = \epsilon |\nabla' W| +\gamma$. Then  one can upgrade $T_1$ to 
\begin{equation}
T_1= \epsilon  \min (\mu\gamma \zeta^{-1},\mu, \rho^{1-\kappa} \mu^K\gamma^K)
\label{2-28}
\end{equation}
Really, $\mu|{\frac {d\ }{dt}}x_1|\le C\zeta $ and one can select direction of time to replace $\rho$ by $\rho^{1-\kappa}$. Then the total contribution to the remainder of all partition elements with $\zeta\le \gamma^\kappa$ is $O(\mu^{-1}h^{-3})$. 

On the other hand, as $\zeta \ge \gamma^\kappa$ and $\gamma\le \mu^{-1+\kappa}$ one can obviously take $T_1=\epsilon \mu \zeta$ and then the total contribution of such partition elements to the remainder does not exceed $C\mu^{-1}\zeta^{-1}h^{-3}\gamma \ll \mu^{-1}h^{-3}$.

\medskip
\noindent
(iv) Finally, let us consider now zone where 
\begin{equation}
 |x_1| \le \gamma=C\mu^{-1}.
 \label{2-29}
%\label{2-38}
\end{equation}
Then I am not getting rid off the cubic terms and should take $T_1=\epsilon \rho$ which one can easily upgrade to $T_1=\epsilon \rho^{1-\kappa}$\,\footnote{\label{foot-7} One can see easily that in this zone ${\frac {d\ }{dt}}b_j=\{a,b_j\}=O(1)$ and $({\frac {d\ }{dt}})^2b_j=\{a,\{a,b_j\}\}=O(1)$; one can prove that the  microhyperbolicity is preserved with respect to the same vector $\ell\rho^{-1}$.} leading to the contribution of this zone to the remainder $C\mu^{-1} h^{-3}$ where factor $\mu^{-1}$ is the measure of zone defined by (\ref{2-29}).
\end{proof}

\subsubsection{}\label{sect-2-1-3} Now let us allow $f_1f_2^{-1}$ to have critical points; however  for the generic magnetic field these critical points are not resonances.

\begin{proposition}\label{prop-2-2} Let conditions $f_1\ne f_2$, $(\ref{2-2})$ and  
\begin{phantomequation}\label{2-30}\end{phantomequation}
\begin{multline}
|\nabla (f_1f_2^{-1})|\le \epsilon_0 \implies\\ \Hess  (f_1f_2^{-1})\, \text{has at least $q$ eigenvalues with absolute values greater than $\epsilon_0$}
\tag*{$(2.30)_q$}\label{2-30-q}
\end{multline}
be fulfilled at $\supp\psi$ with $q\ge3$.  Then under condition $(\ref{2-1})$  the remainder is $O(\mu^{-1}h^{-3})$ while the main part of asymptotics is given by  Weyl formula.
\end{proposition}

\begin{proof} Let $\gamma =\epsilon |\nabla (f_1f_2^{-1})|+{\bar\gamma}$, ${\bar\gamma}=C\mu^{-1}$. Let us consider $\gamma$-admissible partition with respect to $x$. 

Then similarly to part (i) of the proof of proposition \ref{prop-2-1} one can take $T_1=\epsilon \mu \rho ^{2-\delta'}\gamma$ with $\rho$ defined by (\ref{2-8}).  Note that
\begin{equation}
|\alpha\nabla( f_1V^{-1})+(1-\alpha)\nabla(f_2V^{-1})|\ge \epsilon \gamma |\alpha -{\bar\alpha}(x)|\qquad \forall \alpha: 0\le \alpha\le 1 
\label{2-31}
\end{equation}
and therefore  the measure in $\xi$-space gets a factor $\gamma^{-1}$ but the measure in $x$-space gets a factor $\gamma^q$ due to condition condition \ref{2-30-q}. 

Then (\ref{2-11}) and (\ref{2-23}) are replaced by the similar estimate with the right hand expression $C\rho \gamma^{q-1}h^{-3}$.

So, the total contribution $(\rho,\gamma)$-elements to the remainder does not exceed
\begin{equation}
C\mu^{-1}h^{-3}\rho \gamma^{q-1} \times \rho^{\kappa-1}\gamma^{-1}\asymp C\mu^{-1}h^{-3}\gamma^{q-2}\rho^\kappa.
\label{2-32}
\end{equation}
Summation with respect to $(\rho,\gamma)$ results in $C\mu^{-1}h^{-3}$ as 
$q\ge 3$ and in $C\mu^{-1}h^{-3}|\log \mu|$ as $q=2$ and this is the total contribution of zone   $\{\rho \gamma \ge C\mu^{-2}\}$.

On the other hand, contribution of zone $\{\rho \gamma \le C\mu^{-2}\}$ to the remainder is $O(\mu h^{-3}\times \mu^{-2})$ since its measure is $O(\mu^{-2})$ under condition $(\ref{2-30})_3$; under condition $(\ref{2-30})_2$ an extra logarithmic factor appears as well.
\end{proof}

\begin{remark}\label{rem-2-3} Probably one can get rid off logarithmic factors as $q=2$ and derive some estimate as $q=1$. I leave it to the curious reader.
\end{remark}

\subsubsection{}\label{sect-2-1-4} 
Now let us consider the vicinity of $\Sigma=\{x:\ f_1=f_2\}$.  

\begin{proposition}\label{prop-2-4} Let $|f_1-f_2|\asymp \dist (x, \Sigma)$ where $\Sigma$ is a 2-dimensional manifold. Let $\psi$ be supported in the small vicinity of $\Sigma$. Then under condition $(\ref{2-1})$  the remainder is $O(\mu^{-1}h^{-3})$ while the main part of asymptotics is given by  Weyl formula.
\end{proposition}

\begin{proof} Let us introduce a scaling function 
\begin{equation}
\gamma = \epsilon \dist (x,\Sigma)+{\frac 1 2}{\bar\gamma}, \qquad {\bar\ell}=C\mu^{-1}.
\label{2-33}
\end{equation}
Then 
\begin{equation}
|\alpha\nabla( f_1V^{-1})+(1-\alpha)\nabla(f_2V^{-1})|\ge \epsilon  |\alpha -{\bar\alpha}(x)|\qquad \forall \alpha: 0\le \alpha\le 1 
\label{2-34}
\end{equation}
holds with $\gamma$-admissible ${\bar\alpha}$ and  also  (\ref{2-8}),(\ref{2-9})  hold. Therefore  one can take $T_1=\epsilon\rho^{1-\kappa}\gamma \mu $. Then the contribution of all $(\gamma,\rho)$ elements to the remainder does not exceed (\ref{2-32}) with $q=3$\,\footnote{\label{foot-8} In comparison with proposition \ref{prop-2-2} there is no factor $\gamma$ in the right hand expression of (\ref{2-34}) and therefore no factor $\gamma^{-1}$ in the estimate of the measure in $\xi$ space.} and summation over $\{\gamma\ge C\mu^{-1}, \rho \gamma\ge C\mu^{-2}\}$ results in $O(\mu^{-1}h^{-3})$. 

Meanwhile the contributions of zones $\{x: \gamma\ge C\mu^{-1}, \rho \gamma\le C\mu^{-2}\}$, $\{x: \gamma\le C\mu^{-1}\}$ to the remainder are $O(\mu h^{-3}\times \mu^{-2})=O(\mu^{-1}h^{-3})$ since the measures of these zones are $O(\mu^{-2})$.
\end{proof}

\subsubsection{}\label{sect-2-1-5}
Summarizing what is proven one gets

\begin{proposition}\label{prop-2-5} Let $(g^{jk})$ be fixed and then $(V_j)$ be generic, more precisely:

\medskip
\noindent
{\rm (i)} Outside of $\Sigma=\{x:\ f_1=f_2\}$ critical points of $f_1f_2^{-1}$ satisfy $(\ref{2-30})_3$ and $(\ref{2-2})$;

\medskip
\noindent
{\rm (ii)} $\Sigma$ be smooth 2-dimensional manifold and $|f_1-f_2|\asymp \dist(x,\Sigma)$.

Finally, let $V$ be general but satisfying $(\ref{0-8})$ at $\supp\psi$. Then under condition $(\ref{2-1})$  the remainder is $O(\mu^{-1}h^{-3})$ while the main part of asymptotics is given by  Weyl formula.
\end{proposition}

\begin{proposition}\label{prop-2-6} In frames of proposition \ref{prop-2-5} Weyl and Magnetic Weyl expressions differ by (far less than) $O(\mu^{-1}h^{-3})$.
\end{proposition}

\begin{proof} An easy proof is left to the reader.
\end{proof}

\begin{remark}\label{rem-2-7} Definitely $M$ in condition (\ref{2-2}) ``critical points of $f_1f_2^{-1}$ are not resonances of order not exceeding $M$\,'' should not be too large and I leave to the curious reader to investigate it.
\end{remark}

\section{Stronger Magnetic Field Case: Estimates}
\label{sect-3}

\subsection{Canonical form}
\label{sect-3-1}
From now one can assume that
\begin{equation}
h^{-\delta_0}\le \mu \le ch^{-1}
\label{3-1}
\end{equation}
with some small fixed exponent $\delta_0>0$. Then I can reduce operator to a canonical form (depending on additional assumptions) and also make decomposition with respect to Hermite functions,  thus arriving to 2-parametric matrices of 2D $\mu^{-1}h$-PDOs $\cA_{pn}(x',\mu^{-1}hD')$ where ere and below $x'=(x_1,x_2)$.

More precisely, assuming that there are no resonances of order not exceeding $M$:
\begin{equation}
|kf_2-lf_1|\ge \epsilon \qquad \forall (k,l)\in \bZ^{+\,2}: k+l\le  M,
\label{3-2}
\end{equation}
a canonical form contains diagonal elements
\begin{multline}
\cA_{pn}= {\frac 1 2}\Bigl( f_1^\# (2p+1) \mu h + f_2^\# (2n+1)\mu h -V^\# +\\
\sum_{l+m+k+j\ge 2} b_{lmkj} \bigl((2p+1)\mu h\bigr)^l \bigl((2n+1)\mu h\bigr)^m  \mu^{3-2l-2m-2k-j}h^j\Bigr)
\label{3-3}
\end{multline}
with $f_j^\#=f_j^\#(x',\mu^{-1}hD')$, $V^\#=V^\#(x',\mu^{-1}hD')$,
$b_{lmkj}=b_{lmkj}(x',\mu^{-1}hD')$
while all  non-diagonal elements are $O(\mu^{2-M})$.

I will discuss later an alternative form as $M=2$.

\subsection{General estimates at regular points}
\label{sect-3-2}

Assume first that there are no resonances of order not exceeding large $M=M(\delta_0)$. Then under condition (\ref{3-1}) perturbation $O(\mu^{2-M})$ is negligible and (\ref{3-3}) is a true diagonal canonical form (with a negligible perturbation). 

\subsubsection{}
\label{sect-3-2-1} In this case an analysis is easy:

\begin{proposition}
\label{prop-3-1} Let  there be no resonances of order not exceeding (large enough) $M$ and condition $(\ref{1-19})$  be fulfilled at $\supp \psi$. Then 
under condition $(\ref{3-1})$

\medskip
\noindent
{\rm (i)} The standard implicit asymptotic formula $(\ref{0-9})$ holds with the remainder estimate $O(\mu^{-1}h^{-3}+\mu^2h^{-2})$.

\medskip
\noindent
{\rm (ii)} In particular,  remainder estimate is $O(\mu^{-1}h^{-3})$ as
$\mu \le h^{-1/3}$.
\end{proposition}

\begin{proof} (i)  Let us for each \emph{pair\/} $(p,n)$ introduce scaling function 
\begin{equation}
\rho_{pn}=\epsilon \bigl(|\cA_{pn}|+|\nabla\cA_{pn}|^2\bigr)^{1/2}+{\bar\rho},\qquad 
{\bar\rho}=C(\mu^{-1}h|\log h|)^{1/2}+C\mu^{-2}
\label{3-4}
\end{equation}
and a corresponding partition. Then 

\begin{claim}\label{3-5}
The contribution of each \emph{group\/} $(p,n,{\sf element})$ to the main part of asymptotics is  $\lesssim \mu^2h^{-2}\rho^4$\,\footnote{\label{foot-9} And often enough is is of this amplitude, so summation results in the correct magnitude of the main part.}.
\end{claim}
On the other hand, one can take 
\begin{equation}
T_1=\epsilon \mu
\label{3-6}
\end{equation}
since the propagation speed is of magnitude $\mu^{-1}$ and also
\begin{equation}
T_0=Ch\rho^{-2}|\log h|.
\label{3-7}
\end{equation}
Note that $T_0\le \epsilon_0 \mu^{-1}$ provided
\begin{equation}
\rho\ge {\bar\rho}_1= C(\mu h|\log h|)^{1/2}
\label{3-8}
\end{equation}
and in this case one can trade $T_0$ to ${\bar T}=Ch|\log h|$ in (\ref{0-9}) as $\sum Q_\iota=I$; there will be also the correction term arising from elements  failing condition (\ref{3-8}); see section \ref{sect-4}. Moreover, 

\begin{claim}\label{3-9}
In the estimate of expression (\ref{0-10})  one can replace $T_0$ by $T_0^*=Ch\rho^{-2}$. 
\end{claim}
So, the contribution of an element  to the remainder does not exceed
\begin{equation}
C\mu ^2 h^{-2}\rho^4 \times h\rho^{-2} \times \mu^{-1}\times \bigl(\rho^2(\mu h)^{-1}+1\bigr) \times \bigl(\rho (\mu h)^{-1}+1\bigr)\times \rho^{-4}
\label{3-10}
\end{equation}
where  $Ch\rho^{-2}|\log h|$ and $Ch\rho^{-2}$   play the roles of $T_0$ and $T_0^*$ in formulae (\ref{0-9}) and (\ref{0-10}) respectively. Also the numbers of indices $n$\ \footnote{\label{foot-10} For which ellipticity is violated for a given $p$.} and $p$\ \footnote{\label{foot-11} Such that $\rho_{pn}\asymp \rho$ as $n$ violates ellipticity.} are estimated by 
$C\bigl(\rho^2(\mu h)^{-1}+1\bigr)$ and $C\bigl(\rho (\mu h)^{-1}+1\bigr)$ respectively. 

One can rewrite  expression (\ref{3-9})  as 
\begin{equation}
C\mu ^{-1} h^{-3} \rho + Ch^{-2}\rho^{-1}+C\mu  h^{-1} \rho^{-2}; 
\label{3-11}
\end{equation}
then  in the zone
\begin{equation}
\bigl\{\rho \ge {\bar\rho}\Def C(\mu^{-1}h|\log h|)^{1/2}\bigr\}
\label{3-12}
\end{equation}
the first term sums to $C\mu^{-1}h^{-3}$ while the last two terms sum to their values as $\rho={\bar\rho}$, which are   $O(h^{-5/2}\mu^{1/2})=O(\mu^{-1}h^{-3}+\mu^2 h^{-2})$ for sure and $O(\mu^2h^{-2})$ respectively.

\medskip
\noindent
(ii) Meanwhile with the same main part the total contribution to the remainder of all groups with $\{\rho \le{\bar\rho}\}$  trivially does not exceed
\begin{equation}
C\mu ^2 h^{-2} \bigl({\bar\rho}^2(\mu h)^{-1}+1\bigr)\times \bigl({\bar\rho}(\mu h)^{-1}+1\bigr)\le C\mu h^{-3}{\bar\rho} + C\mu^2 h^{-2}.
\label{3-13}
\end{equation}
Obviously, this  expression is $O(\mu^{-1}h^{-3})$ as $\mu \le C(h|\log h|)^{-1/3}$ and $O(\mu^2h^{-2})$ as $\mu \ge Ch^{-1/3}|\log h|^{1/3}$.

\medskip
\noindent
(iii) To finish the proof I need to reconsider contribution to the remainder of the elements with $\{\rho_{pn}\le {\bar\rho}\}$ in the border case
\begin{equation}
h^{-1/3}|\log h|^{-1/3}\le \mu \le h^{-1/3}|\log h|^{1/3}.
\label{3-14}
\end{equation}
Let us introduce another scaling function 
\begin{equation}
\varrho= \epsilon |\nabla^2 \cA_{pn}|+{\frac 1 2}{\bar\varrho}, \qquad {\bar\varrho}=|\log h|^{-K},
\label{3-15}
\end{equation}
calculated with $(p,n)$ delivering minimum to  $\rho_{pn}$ and let us introduce the corresponding partition. 

Then for any element with ${\varrho}\ge {\bar\varrho}$ one can calculate easily that the relative measure of the zone $\{(x',\xi'):\ \min_{p,n}\rho_{pn} \le C{\bar\rho}\}$ is $O({\bar\rho}|\log h|^K)$ and then the total contribution of zone $\{(x',\xi'):\ \rho \le {\bar\rho},\ \varrho \ge C{\bar\varrho}\}$ to the remainder is much less than $O(\mu^{-1}h^{-3})$. 

On the other hand, let us consider elements with $\varrho\le C{\bar\varrho}$.
Since the total contribution of subelements  with $\rho\le \rho^*_1\Def C(\mu h)^{1/2}$ to the remainder is estimated properly, one needs to consider only subelements with $\rho^*_1  \le \rho\le {\bar\rho}_1$.

But on such subelements $\rho\varrho^{-1}$ is a scaling function as well and using it one can easily decrease $T_0$ to  $Ch|\log h|{\bar\varrho}\rho^{-2}$ leaving $T_1=\epsilon \mu$ ; this will add an extra factor $|\log h|^{2-K}$ to the estimate of the contribution of this zone to the remainder and this factor leads to the needed estimate $O(\mu^{-1}h^{-3})$.
\end{proof}

\subsubsection{}
\label{sect-3-2-2} Assume now that the is a critical point of $f_1f_2^{-1}$:

\begin{proposition}
\label{prop-3-2} Let  there be no resonances of order not exceeding (large enough) $M$ and condition $(\ref{2-30})_3$ be fulfilled at $\supp \psi$. Then 
under condition $(\ref{3-1})$

\medskip
\noindent
{\rm (i)} The standard implicit asymptotic formula holds with the remainder estimate $O(\mu^{-1}h^{-3}+\mu^2h^{-2})$.

\medskip
\noindent
{\rm (ii)} In particular, the remainder estimate is $O(\mu^{-1}h^{-3})$ as $\mu\le C h^{-1/3}$.
\end{proposition}

\begin{proof} (i) The arguments of the proof of proposition \ref{prop-3-1} still work without condition (\ref{1-19}) with the exception of the estimate $C\bigl(\rho (\mu h)^{-1}+1\bigr)$ of the number of the indices ``$p$''. 

However, let us introduce another scaling function
\begin{equation}
\gamma \Def \epsilon_1 |\nabla (f_1f_2^{-1})|
\label{3-16}
\end{equation}
and if on some group $\rho_{pn}\le \gamma$ then the number of indices ``$p$'' should be estimated by $C\bigl(\rho (\gamma \mu h)^{-1}+1\bigr)$; otherwise this number should be estimated by $C\bigl( ( \mu h)^{-1}+1\bigr)$. 

Note that if condition $\rho_{pn}\lesssim \gamma$ (or $\rho_{pn}\gtrsim \gamma$) is fulfilled at some point of $\rho$-element, then the same condition (with another implicit constant) is fulfilled at any other point of this element.

Anyway, this modification adds no more than one factor $\gamma^{-1}$ to the estimate, but the factor $\gamma^q$ comes from condition \ref{2-30-q}; then summation over $\gamma$-partition results in the same estimate $O(\mu^{-1}h^{-3}+\mu^2h^{-2})$ as $q\ge 2$ and in the same estimate but with an extra logarithmic factor as $q=1$. Since I do not use  microlocal $\gamma$-partition, no estimate of $\gamma $ from below  is needed.

Therefore,
\begin{claim}\label{3-17}
The total contribution of elements with $\rho_{pn}\ge {\bar\rho}$ to the remainder is $O(\mu^{-1}h^{-3}+\mu^2h^{-2})$.
\end{claim}

\medskip
\noindent
(ii) Meanwhile, the contribution of zone $\{x:\ \rho \le{\bar\rho},\ |\nabla (f_1f_2^{-1})| \le \gamma\}$ to the remainder does not exceed 
$C\mu h^{-3}{\bar\rho}\gamma^{q-1}+C\mu^2h^{-2}\gamma^q$ which sums to $O(\mu^{-1}h^{-3}+\mu^2h^{-2})$ 
except in the border case (\ref{3-14}).  However even in  this case $q=2$ only a contribution of zone $\{x: |\nabla (f_1f_2^{-1})|\ge |\log h|^{-1}\}$ 
is not properly estimated. However one can tackle it by the same arguments as in the part (iii) of the proof of proposition \ref{prop-3-1}.
\end{proof}

\begin{remark}\label{rem-3-3} Probably one can get rid off logarithmic factors as $q=1$.
\end{remark}

\subsection{Sharp asymptotics at regular points}
\label{sect-3-3}

The purpose of this and the next odd-numbered subsections is to consider the case
\begin{equation}
ch^{-1/3}\le \mu \le ch^{-1}
%\label{3C-44}
%\label{3A-49}
\label{3-18}
\end{equation}
and derive some under non-degeneracy condition remainder estimate $O(\mu^{-1}h^{-3})$. 

\subsubsection{}
\label{sect-3-3-1} Note first that

\begin{claim}\label{3-19}%\label{3C-45}%\label{3C-58}
In the generic case  $\nabla (f_1f_2^{-1})({\bar x})=0$ implies $\nabla (Vf_1^{-1})({\bar x})\ne 0$ and then microhyperbolicity condition holds with any $\ell$ such that $\langle \ell, \nabla (Vf_1^{-1})\rangle>0$
\end{claim}
and therefore the following statement is generic:

\begin{claim}\label{3-20}%\label{3C-46}
Let $\nabla (f_1f_2^{-1})({\bar x})=0$ and $\nabla (Vf_1^{-1})({\bar x})\ne 0$. Then for any $\psi$ supported in the small vicinity of ${\bar x}$ asymptotics with the magnetic Weyl main part and remainder estimate $O(\mu^{-1}h^{-3})$ holds.
\end{claim}

\subsubsection{}
\label{sect-3-3-2} Therefore in what follows one can assume that $\nabla (f_1f_2^{-1})$ is disjoint from 0. Further, one should consider only vicinities of points  where microhyperbolicity condition (\ref{1-15}) is violated. 

\begin{proposition}\label{prop-3-4}%\label{prop-3C-6}
Let conditions $(\ref{0-8})$ and $(\ref{3-18})$ be fulfilled. Moreover, let us assume that on $\supp\psi$ there are no resonances of order not exceeding $M$ and also $|\nabla (f_1f_2^{-1})|\ge \epsilon_0$. Furthermore, let condition
\begin{phantomequation}\label{3-21}\end{phantomequation}
\begin{multline}
\nu (\rho)\Def \mes \bigl\{ (x,\alpha):\ 0\le \alpha\le 1,\\|\nabla\bigl( \alpha\log f_1+(1-\alpha)\log f_2- \log V\bigr)|\le \rho\bigr\}=O(\rho^q)\qquad
\text{as\ }\ \rho\to+0
\tag*{$(3.21)_q$}\label{3-21-q}%\label{3-18-q}
\end{multline}
be fulfilled with $q>1$. Then  asymptotics with the standard implicit main part $(\ref{0-9})$ and the remainder estimate 
\begin{phantomequation}\label{3-22}\end{phantomequation}
\begin{equation}
O\bigl(\mu^{-1}h^{-3}+\mu^2 h^{-2}(\mu^{-1}h)^{(q-1)/2}|\log h|^{(q+1)/2}\bigr)
%\label{3C-62}
\tag*{$(3.22)_q$}\label{3-22-q}\end{equation}
holds.
\end{proposition}

\begin{proof}
 As $T=Ch\rho^{-2}|\log h|$ the contribution to the remainder of $\rho$-elements does not exceed (\ref{3-10}) multiplied by $\nu (\rho)\rho^{-1}$
\begin{multline*}
C\mu  h^{-1} \rho^{-2} \bigl(\rho^2(\mu h)^{-1}+1\bigr) \times \bigl(\rho (\mu h)^{-1}+1\bigr)\times \nu(\rho) \rho^{-1}\asymp \\
C\mu h^{-1}\bigl( \rho^3(\mu h)^{-2}+\rho (\mu h)^{-1} +1\bigr) \nu(\rho)\rho^{-3}\le\\
C\mu^{-1} h^{-3}  \rho^q + Ch^{-2} \rho^{q-2} +C \mu h^{-1}\rho^{q-3}
\end{multline*}
where the last factor in the left-hand expression is the total measure of $\rho$-elements.

Here the first term in the right-hand expression always sums to $O(\mu^{-1}h^{-3})$ while the second and the third factor sum to their values as $\rho=1$ (i.e. $O(\mu^{-1}h^{-3})$ for sure) plus their values as $\rho={\bar\rho}=C(\mu^{-1}h|\log h|)^{1/2}$ which do not exceed the second term in \ref{3-22-q}. As $q=2,3$ the second or the third term respectively acquires an extra logarithmic factor  but it does not change the estimate.

Furthermore, contribution of zone $\bigl\{\rho \le {\bar\rho}\bigr\}$ to the remainder does not exceed $C\mu h^{-3}\nu ({\bar\rho})+ C\mu^2h^{-2}\nu ({\bar\rho}){\bar\rho}^{-1}$ which is the second term in \ref{3-22-q}.
\end{proof}

\begin{corollary}\label{cor-3-5}  In frames of proposition \ref{prop-3-4}  remainder estimate $O(\mu^{-1}h^{-3})$ holds as   $q>3$.
\end{corollary}

\begin{remark}\label{rem-3-6}
In the generic case condition $(\ref{3-21})_4$ is fulfilled.
\end{remark}

\subsubsection{}
\label{}
Let us consider a special case 
\begin{equation}
\epsilon h^{-1}\le \mu \le ch^{-1}.
%\label{3C-63}
\label{3-23}
\end{equation}

\begin{proposition}\label{prop-3-7}
Let conditions $(\ref{0-8})$ and $(\ref{3-23})$ be fulfilled. Moreover, let us assume that on $\supp\psi$ there are no resonances of order not exceeding $M$ and also  $|\nabla (f_1f_2^{-1})|\ge \epsilon_0$. Furthermore, let condition
\begin{phantomequation}\label{3-24}\end{phantomequation}
\begin{multline}
|(2p+1)\mu h f_1V^{-1}+(2n+1)\mu h f_2V^{-1}-1|+\\
|\nabla\bigl( (2p+1)\mu h f_1V^{-1}+(2n+1)\mu h f_2V^{-1}-1\bigr)|\le \epsilon_0 \implies\\ \Hess  \bigl((2p+1)\mu h f_1V^{-1}+(2n+1)\mu h f_2V^{-1}-1\bigr)  \text{\ \ has at least $r$ eigenvalues  }\\ \text{with absolute values greater than $\epsilon_0$}
\tag*{$(3.24)_r$}\label{3-24-r}
\end{multline}
be fulfilled with $r\ge 1$. Then asymptotics with the standard implicit main part $(\ref{0-9})$ and the remainder estimate $(\ref{3-22})_{r+1}$  i.e.
$O(h^{-2}+h^{r-4}|\log h|^{(r+2)/2})$ holds.
\end{proposition}

\begin{proof} Follows easily the proof of proposition \ref{prop-3-4}.
\end{proof}

\begin{remark}\label{rem-3-8}
One can get rid off the logarithmic factors in estimate \ref{3-22-q} but I do not care since I am interested only in the generic cases $q=4$, $r=4$.
\end{remark}

\subsection{General asymptotics at resonances}
\label{sect-3-4}

Now assume that there are $(k,l)$ resonances of order $m=k+l\ge 3$. However considering $\epsilon_1$-vicinity of any point one can assume that

 \begin{claim}\label{3-25}%\label{3C-15}
There are no $(k',l')$-resonances with $k'+l'\le M$ unless $k'/k = l'/l \in \bZ$.
\end{claim}

\subsubsection{}
\label{sect-3-4-1} So, let us consider  
\begin{equation}
\Xi_{kl}=\{x:\ f_1k=f_2l\}
\label{3-26}
\end{equation}
which under condition (\ref{3-4}) is a smooth surface. One can assume without any loss of the generality that the analogue of (\ref{2-20}) holds: $kf_1-lf_2=x_1$ while $(x_2,\xi_1,\xi_2)$ are coordinates on $\Xi_{kl}$. Let us introduce a scaling function 
\begin{equation}
\gamma =\epsilon |x_1|+{\frac 1 2}{\bar\gamma},\qquad
{\bar\gamma}=\mu^{-1+\delta}
\label{3-27}
\end{equation}
with arbitrarily small exponent $\delta>0$.

\begin{proposition}\label{prop-3-9} Let conditions $(\ref{3-25})$ with $k+l=m\ge3$   and $(\ref{1-19})$ be fulfilled at $\supp\psi$. 
Then under condition $(\ref{3-1})$ asymptotics with the standard implicit main part $(\ref{0-9})$ and the remainder estimate  $O(\mu^{-1}h^{-3}+\mu^2h^{-2})$ holds.
\end{proposition}

\begin{proof}[Proof, Part I]  I will give the proof working the worst-case scenario $m=3$; as $m\ge 4$ one can simplify the proof. In this part I am going to prove that
\begin{claim}\label{3-28}
The contribution of zone $\{|x_1|\ge {\bar\gamma}\}$ to the remainder is $O(\mu^{-1}h^{-3}+\mu^2h^{-2})$.
\end{claim}
After reduction to precanonical form with non-diagonal terms corresponding to resonances, in this zone one can get rid off non-diagonal terms (modulo $O(\mu^{-M})$).  Let us consider scaling function $\rho_{pn}$ introduced by (\ref{3-4}) for a full symbol of $\cA_{pn}$.

\medskip
\noindent 
(i) Consider first subzone
\begin{equation}
\bigl\{x:\ |x_1|\ge \max(\rho, {\bar\gamma})\bigr\}.
\label{3-29}
\end{equation}
Then one can apply the same arguments as in the proof of proposition \ref{prop-3-1}; however there is a problem\footnote{\label{foot-12} It is not a problem at all as $m=5$ and rather a marginal problem as $m=4$ but for $d=3$ this requires a certain attention.}:   as function of $x_1$ \ $\rho_{pn}$ remains $\gamma$-admissible only as 
\begin{equation}
\gamma\ge {\bar\gamma}_1\Def \mu^{(4-2m)/3}.
\label{3-30}
\end{equation}
So far the arguments as in the proof of proposition \ref{prop-3-1} result in\footnote{\label{foot-13} The border case (\ref{3-14}) (may be with the different powers of logarithm) should be covered only in zone $\{|\log h|^{-K_1}\le |x_1|\le \epsilon\}$.}

\begin{claim}\label{3-31}
The contribution of zone $\bigl\{x:\ |x_1| \ge \max(\rho, {\bar\gamma}_1)\bigr\}$ with ${\bar\gamma}_1= \mu^{-2/3}$ to the remainder is $O(\mu^{-1}h^{-3}+\mu^2h^{-2})$.
\end{claim}

\medskip
\noindent
(ii) Consider now subzone
\begin{equation}
\bigl\{\gamma \ge {\bar\gamma},\ \rho\ge \gamma ,\bigr\}.
\label{3-32}
\end{equation}
Note that here the derivatives of $\cA_{pn}$ with respect to $x_1$ and to $(x_2,\xi_1,\xi_2)$ have different values: while the derivative with respect to $x_1$ measures speed with respect to $\xi_1$ and the shift with respect to $\xi_1$ is quantum observable as 
$|\nabla_{x_1}\cA_{pn}|\times \gamma \ge C\mu^{-1}h|\log h|$  provided $\gamma$ is the scale with respect to $x_1$, other derivatives measure speed with respect to $(\xi_2, x_1,x_2)$ and the shift is quantum observable as 
$|\nabla'\cA_{pn}|\times \zeta \ge C\mu^{-1}h|\log h|$,
provided $\zeta$ is the scale with respect to these variable where here and below
$\nabla'\Def \nabla_{x_2,\xi_1,\xi_2}$.

So, let us introduce the third scaling function 
\begin{equation}
\zeta_{pn}=\epsilon \bigl(|\cA_{pn}|+|\nabla'\cA_{pn}|^2\bigr)^{1/2}+
(\rho\gamma)^{1/2}
\label{3-33}
\end{equation}
where the last term will be actually included later.

Then one can make a $\zeta$-admissible partition. So I have now $(\gamma,\rho,\zeta)$ elements with $\gamma\le \zeta\le\rho$ where $\zeta$ is the scale with respect to $(x_2,\xi_1,\xi_2)$ and $\gamma$ is the scale with respect to $x_1$ while $\rho$ at this moment lost its scaling role. However while $\zeta$ is $\gamma$-admissible function with respect to $x_1$, $\rho$ preserves its magnitude in as $x_1$ varies by $O(\gamma)$ only under condition
\begin{equation}
\rho \ge \mu^{-2}\gamma^{-2}.
\label{3-34}
\end{equation}
I claim that  that then 
\begin{align}
&T_0=Ch|\log h| \bigl(\rho\gamma +\zeta^2\bigr)^{-1},\qquad
T_0^*=Ch \bigl(\rho\gamma +\zeta^2\bigr)^{-1}\label{3-35}\\
\intertext{while}
&T_1= \epsilon\mu  \min\bigl( \zeta\rho^{-1},\gamma\zeta^{-1}\bigr)\asymp \epsilon\mu \zeta\gamma\bigl( \rho\gamma+\zeta^2\bigr)^{-1}.\label{3-36}
\end{align}
Really, as  $\rho \gamma \asymp \zeta^2$ propagation speed with respect to  $\xi_1$ is $\asymp\rho$ it is dual to $x_1$ of the scale $\gamma$.  Meanwhile speeds with respect to all other variables are bounded by $C\zeta$ and for given $T_1$ magnitudes of $\zeta,\gamma,\rho$ are preserved. 

On the other hand, as $\rho \gamma \le \epsilon_1\zeta^2$ propagation speed with respect to $(x_1,x_2,\xi_2)$ is $\asymp \zeta$ and they are dual to $(\xi_1,\xi_2,x_2)$ of the scale $\zeta$ while propagation speed with respect to $\xi_1$ is bounded by $C\rho$ and for given $T_1$ magnitudes of $\zeta,\gamma,\rho$ are preserved. 

I leave to the reader the standard justification on the quantum level (based on energy estimates approach).

Thus $T_0^*T_1^{-1}\asymp \mu ^{-1}h \gamma^{-1}\zeta^{-1}$; in virtue of the last term in the definition of $\zeta$   one can skip $\rho\gamma$ in $(\zeta^2+\rho\gamma)$ here and below.

Then the  total contribution of all $(\rho,\zeta,\gamma)$ elements to the remainder does not exceed
\begin{equation}
C\mu^2h^{-2}\times \gamma \times T_0^*T_1^{-1}\times
\bigl((\rho\gamma+\zeta^2)(\mu h)^{-1}+1\bigr)\times \bigl(\rho (\mu h)^{-1}+1\bigr)
\label{3-37}
\end{equation}
where $C\bigl((\rho\gamma+\zeta^2)(\mu h)^{-1}+1\bigr)\asymp C\bigl((\zeta^2(\mu h)^{-1}+1\bigr)$ is an upper bound for a number of indices $n$ violating ellipticity for a given index $p$. 

If one picks up only ``$1$'' from both factors with the parentheses in (\ref{3-37}) and replaces $T_0T_1^{-1}$ by $1$, then summation with respect to partitions results in $C\mu^2h^{-2}$; on the other hand, since $\zeta^2\le \rho$, one can rewrite the above expression (\ref{3-37}) as
\begin{equation}
C\mu h^{-1}\times \zeta^{-1}
\bigl(\zeta^2 (\mu h)^{-1}+1\bigr)\times \rho (\mu h)^{-1}\asymp
C h^{-2}  \bigl(\zeta  (\mu h)^{-1}+\zeta^{-1}\bigr)\rho .
\label{3-38}
\end{equation}
Then summation with respect to $\zeta$ from $(\rho \gamma)^{1/2}$ to $\rho$ results in 
\begin{equation}
Ch^{-2}\rho 
\bigl(\rho  (\mu h)^{-1}+\rho^{-1/2}\gamma^{-1/2}\bigr) \asymp 
C\mu^{-1} h^{-3}\rho^2 + 
Ch^{-2}\rho^{1/2}\gamma^{-1/2}.
\label{3-39}
\end{equation}
The second term in the right hand expression sums with respect to $(\rho,\gamma)$ to $Ch^{-2}{\bar\gamma}^{-1/2}\ll C\mu^2 h^{-2}$. However the first term sums to $C\mu^{-1}h^{-3}|\log h|$ and this logarithmic factor appears due to summation with respect to $\gamma$.

To get rid off this factor let us notice that only case $\mu \le h ^{-1/3}|\log h|^{1/3}$ needs to be addressed and only zone $\{\zeta \ge \gamma^\kappa\}$ should be reconsidered (with an arbitrarily small exponent $\kappa>0$); in this case $\zeta \gg (\rho\gamma)^{1/2}$. Moreover,  only term $Ch^{-4}\rho\zeta^2 \gamma\times T_0^*T_1^{-1}$ in (\ref{3-37}) should be reexamined.

However then one does not need to use the canonical form but rather a weak magnetic field approach and  take $T_0^*=h\zeta^{-2}$ and $T_1=\mu \zeta$ and the contribution of this zone to the term in question does not exceed $C\mu^{-1}h^{-3}\int \zeta^{-1}\,d\gamma \le C\mu^{-1}h^{-3}$. Therefore
\begin{claim}\label{3-40}
The contribution of zone 
$\bigl\{ \rho \ge \max (\gamma, \mu^{-2}\gamma^{-2}),\ \gamma \ge {\bar\gamma}\bigr\}$ to the remainder is $O(\mu^{-1}h^{-3}+\mu^2h^{-2})$.
\end{claim}

\medskip
\noindent
(iii) Now let us consider the remaining part of the zone 
$\{\epsilon\ge \gamma \ge {\bar\gamma}\}$. In this zone let us introduce the scaling function 
\begin{equation}
\eta = \mu ^2\gamma^3\rho.
\label{3-41}
\end{equation}
Then $\rho$  as a function of $x_1$ is $\eta$-admissible. Let us modify definition (\ref{3-33}) etc replacing $\gamma$ by $\eta$:
\begin{equation}
\zeta_{pn}=\epsilon \bigl(|\cA_{pn}|+|\nabla'\cA_{pn}|^2\bigr)^{1/2}+
(\rho\gamma)^{1/2}
\tag*{$(3.33)^*$}
\label{3-33-*}
\end{equation}
and then (\ref{3-35}), (\ref{3-36})  are also modified in the same way (anyway, terms $\rho\gamma$ originally and $\rho\eta$ now are not important): 
\begin{equation}
T_1= \epsilon\mu  \min\bigl( \zeta\rho^{-1},\eta\zeta^{-1}\bigr)\asymp \epsilon\mu \eta\zeta^{-1}
\tag*{$(3.35)^*$}
\label{3-35-*}
\end{equation}
and $T_0^*T_1^{-1}\asymp \mu^{-1}h \eta^{-1}\zeta^{-1}$.

Then modified (\ref{3-37}) and (\ref{3-38}) expressions
\begin{equation}
C\mu h^{-1}\times \zeta^{-1}
\bigl(\zeta^2 (\mu h)^{-1}+1\bigr)\times \rho (\mu h)^{-1}\gamma\eta^{-1}\asymp
C h^{-2} \times 
\bigl(\zeta  (\mu h)^{-1}+\zeta^{-1}\bigr)\rho \gamma\eta^{-1}
\tag*{$(3.37)^*$}
\label{3-37-*}
\end{equation}
estimate contribution of all $(\rho,\gamma,\eta,\zeta)$ elements. Here an (unpleasant) factor $\gamma\eta^{-1}$ appears since the measure of the strip remains $\gamma$ and  $\gamma^{-1}$ is replaced by $\eta^{-1}$ in $T_0^*T_1^{-1}$. 

Plugging $\eta$ into \ref{3-37-*} one gets 
\begin{equation}
C\mu^{-3}h^{-3}\gamma^{-2}\zeta +C\mu^{-2}h^{-2} \gamma^{-2}\zeta^{-1}
\label{3-42}
\end{equation}
and the first term sums to $o(\mu^{-1}h^{-3})$.

The second term sums with respect to $\zeta$ to its value at the smallest $\zeta$ satisfying conditions $T_1=\mu \eta \zeta^{-1}\ge T_0=h\zeta^{-2}$ and $\zeta \ge (\rho\eta)^{1/2}$:
\begin{equation}
\zeta = \max \bigl(\mu^{-1}h\eta^{-1}, (\rho\eta)^{1/2}\bigr)= \max \bigl(\mu^{-3}h\rho^{-1}\gamma^{-3}, \rho \mu \gamma^{3/2}\bigr)
\label{3-43}
\end{equation}
so one gets 
\begin{equation*}
C\mu^{-2}h^{-2} \gamma^{-2}\min \bigl(  h^{-1}\rho \mu^3\gamma^3 , \rho^{-1}\mu^{-1}\gamma^{-3/2}\bigr)
\end{equation*}
which sums with respect to $\rho$ to $C\mu^{-1}h^{-5/2} \gamma^{-5/4}$ and then with respect to $\gamma$ to $O(\mu^{1/4} h^{-5/2})\ll \bigl(\mu^{-1}h^{-3}+ \mu^2 h^{-2}\bigr)$.

Meanwhile if $T_1\le T_0$ then one just replaces $T_1$ by $T_0$ and gets 
$C h^{-4}\zeta^2 \rho\gamma + C\mu h^{-3}\rho\gamma$ with $\zeta=(\rho\eta)^{1/2}$ i.e.
\begin{equation}
C h^{-4}(\rho \eta)^{1/2}\rho\gamma  + C\mu h^{-3}\rho\gamma\asymp 
C h^{-4}\mu \rho^2\gamma^{5/2}  + C\mu h^{-3}\rho\gamma
\label{3-44}
\end{equation}
and \emph{only} if $(\rho\eta)^{1/2}\le \mu^{-1}h\eta^{-1}$ or equivalently
$\rho \le  {\bar\rho}\Def \mu^{-2}h^{1/2}\gamma^{-9/4}$. Plugging ${\bar\rho}$ in (\ref{3-44}) I get 
$C \bigl(h^{-3}\mu ^{-3}\gamma^{-2}  + \mu ^{-1}h^{-5/2}\gamma^{-5/4}\bigr)$; summation with respect to $\gamma$ results in $o(\mu^{-1}h^{-3}+\mu^{1/4}h^{-5/2})$. Therefore

\emph{Statement $(\ref{3-28})$ is proven.} 
\end{proof}

\begin{proof}[Proof of Proposition 3.9, Part II]\label{pf-3-9-II}  I am left with zone $\{|x_1|\le {\bar\gamma}\}$; its contribution to the remainder does not exceed $C\mu h^{-3}{\bar\gamma}=C\mu^\delta h^{-3}$ which is $O(\mu^2h^{-2})$ as $\mu\ge h^{-1/2-\delta'}$ and therefore only case $\mu \le h^{-1/2-\delta'}$ needs to be addressed. In this zone I use precanonical form with non-diagonal matrix elements but without singular terms. 

\begin{claim}\label{3-45} 
Let $\cA^0_{pn}$ denote diagonal matrix elements,
\end{claim}
One can see easily that 
\begin{equation}
|\nabla '\cA^0_{pn}|\equiv  |\nabla ' \log (f_jV^{-1})|\quad \mod O(\gamma)\qquad \text{as\ \ } |x_1|\le \gamma
\label{3-46}
\end{equation}
(as $\cA_{pn}$ is non-elliptic) does not actually depend on $p$; let us define
\begin{equation}
\zeta \Def |\nabla ' \log (f_jV^{-1}|\bigr|_{x_1=0} + (\rho \gamma)^{1/2}.
\label{3-47}
\end{equation}
\begin{claim}\label{3-48} 
Let us define $\rho_{pn}$  as before but with $\cA^0_{pn}$ instead of $\cA_{pn}$.
\end{claim}
in contrast to $\zeta$ $\rho_{pn}$ strongly depends on $p$.

Then, as before $\mu^{-1}\rho$ controls the propagation speed with respect to $x_1$ and thus bounds the propagation speed with respect to $\zeta$ while $\mu^{-1}\zeta$ controls the propagation speed with respect to $(x_1,x_2,\xi_2)$ as long as $\zeta\ge C\gamma$, $\rho\ge C\gamma$. I remind that $\gamma\ge C\mu^{-1}$.

However the propagation speed with respect to $\rho$ is a different matter. Considering commutator of 
\begin{equation*}
\cA= f_1V^{-1} Z_1^*Z_1+ f_2 Z_2^*Z_2+2\mu^{-1}\Re \bigl(\beta Z_1^*Z_2^2\bigr) +\dots
\end{equation*}
(assuming that resonance is $2f_2=f_1$) with 
\begin{equation}
\partial_{x_1}\cA^0 = \varkappa Z_1^*Z_1 + \omega \cA^0+\dots, \qquad \varkappa= f_1  \bigl(\partial_{x_1}(\log (f_1f_2^{-1})\bigr)
\label{3-49}
\end{equation}
(with $\dots =O(\gamma)$) one can see easily that 
\begin{equation}
\varkappa \mu^{-1}
\bigl[\cA, \partial_{x_1}\cA^0\bigr] \equiv \bigl[\Re \beta Z_1^*Z_2^2,Z_1^*Z_1\bigr]= 2\varkappa \Re \beta Z_1^*Z_2^2 \quad\mod O(\gamma)
\label{3-50}
\end{equation}
which is bounded by $1$. Therefore 
\begin{claim}\label{3-51}
The propagation speed with respect to $\rho$ does not exceed $1$. Furthermore, $\rho$ is properly defined as $\rho\ge C\mu h|\log h|$  (the logarithmic uncertainty principle).
\end{claim}
Without nondegeneracy condition there is not much use of $\zeta$; let us consider elements with $\rho\ge C\gamma$. Due to (\ref{3-51}) I can pick up $T_1=\epsilon \rho$. Then the contribution to the remainder of all such elements does not exceed 
\begin{equation}
C\mu^2h^{-2}\gamma \times h\zeta^{-2}\rho^{-1}
\bigl(\zeta^2  (\mu h)^{-1} +1\bigr) \times \bigl(\rho  (\mu h)^{-1} +1\bigr)
\label{3-52}
\end{equation}
and as long I include $O(\mu^2h^{-2})$ into final remainder estimate I can skip ``$+1$'' in the last factor (due to the same arguments as before) resulting in
\begin{equation}
C\mu h^{-2}\gamma \zeta^{-2}
\bigl(\zeta^2  (\mu h)^{-1} +1\bigr) \asymp Ch^{-3}\gamma + C\mu h^{-2}\zeta^{-2}\gamma\le Ch^{-3}\gamma + C\mu h^{-2}\rho^{-1};
\label{3-53}
\end{equation}
as $\gamma={\bar\gamma}$ this expression does not exceed $Ch^{-3}{\bar\gamma}+C\mu h^{-2}\rho^{-1}$ which sums with respect to $\rho $ to
$Ch^{-3}{\bar\gamma}|\log h|+ C\mu h^{-2}{\bar\rho}^{-1}$ and the last term is $O(\mu^2h^{-2})$. 

On the other hand, in the zone $\{\rho \le {\bar\rho}\}$ I pick up $T^*_0T_1^{-1}=1$ and its contribution to the remainder does not exceed $C{\bar\rho}^2{\bar\gamma}^2h^{-4}+ C\mu h^{-3}{\bar\rho}{\bar\gamma}+C\mu^2h^{-2}$.

Therefore, as $h^{-1/3-\delta'}\le \mu \le h^{-1/2-\delta'}$,  the contribution to the remainder of the zone $\{|x_1|\le {\bar\gamma}=\mu^{-1+\delta}\}$ with sufficiently small $\delta=\delta(\delta')>0$ does not exceed $C\mu^2h^{-2}$. In this case proposition \ref{prop-3-9} is also proven.
\end{proof}

I will need the following 

\begin{proposition}\label{prop-3-10} Let conditions $(\ref{0-8})$ and $(\ref{3-1})$ and $|f_1-f_2|\ge \epsilon$ be fulfilled. Let us consider the precanonical form. Finally, let $\cQ$ be a $\rho$-admissible partition element element in $(q_1,q_2)$, then quantized as 
$\cQ \bigl((h^2D_3^2+\mu^2x_3^2), (h^2D_4^2+\mu^2x_4^2)\bigr)$ with
\begin{equation}
\rho\ge C\mu h|\log h|+ C{\bar\gamma}_0^k\gamma^{1-k}
\label{3-54}
\end{equation}
and $\psi$ be $\gamma$-admissible with respect to $x_1$, either supported in $\{|x_1|\asymp \gamma\}$ as $\gamma> {\bar\gamma}_0\Def C\mu^{-1}$ or supported in $\{|x_1|\lesssim \gamma\}$ as $\gamma={\bar\gamma}_0$.
Then 
\begin{align}
&|F_{t\to h^{-1}\tau} \chi_T(t)\Gamma \bigr(\psi \cQ u Q ^t_y\bigl)|\le Ch^s,
\label{3-55}\\
\intertext{and}
&|F_{t\to h^{-1}\tau} {\bar\chi}_T(t)\Gamma \bigr(\psi \cQ u Q ^t_y\bigl)|\le  C\rho \gamma h^{-3}
\label{3-56}
\end{align}
as $|\tau|\le \epsilon$
\begin{equation}
Ch|\log h|\rho^{-1}\le T\le \epsilon\mu^{-1}.
\label{3-57}
\end{equation}
\end{proposition}

\begin{proof}
Consider the propagation with respect to either $(x_3,\mu^{-1} hD_3)$ or
$( x_4,\mu^{-1}hD_4)$. Due to (\ref{0-8}) on energy levels close to 0 at least one of $\mu^2x_j^2+h^2D_3^j$ is of magnitude 1 ($j=3,4$). The propagation speed with respect to $(x_3, x_4,\mu^{-1}hD_3,\mu^{-1}hD_4)$ is $\asymp 1$. Therefore one can trade $T\le \epsilon \mu^{-1}$ to ${\bar T}=Ch|\log h|$ and in the estimate of Fourier transform to ${\bar T}^*=Ch$. The remaining part of the proof is easy and left to the reader.
\end{proof}

\begin{proof}[Proof of Proposition 3.9, Part III]\label{pf-3-9-III} Therefore only case $h^{-\delta_0}\le \mu \le h^{-1/3-\delta'}$ remains to be addressed where $\delta_0>0$ is small and fixed and $\delta'>0$ is arbitrarily small. It follows from Parts I,II that exponents $\delta>0$ in the definition of ${\bar\gamma}$ and $\delta'>0$ are \emph{independently\/} small. 

\medskip
\noindent
(i) Let us consider elements with $|\nabla 'V^{-1}f_j|\asymp \varsigma\ge C\gamma$. Then as $T_0= Ch|\log h|(\varsigma^2+\rho\gamma)^{-1}\le \epsilon \rho$, i.e. as
\begin{equation}
\rho \ge \varrho\Def C\min \bigl(h|\log h|\varsigma^{-2}, (\gamma^{-1}h|\log h|)^{1/2}\bigr)+C\gamma,
\label{3-58}
\end{equation}
$\rho $ is preserved on the time interval $T_0$ which can be traded to $T_1=\epsilon\mu \varsigma$. Therefore the contribution to the remainder of all such elements does not exceed
\begin{equation}
C\rho\gamma (\varsigma^2+\mu h) h^{-4}T_0^* T_1^{-1} \asymp C\mu^{-1}h^{-3}\rho\gamma\varsigma^{-1}+ C h^{-2}\rho\gamma\varsigma^{-3};
\label{3-59}
\end{equation}
summation with respect to $\varsigma,\rho, \gamma$ trivially results in $O(\mu^{-1}h^{-3}+\mu^2h^{-2})$.

On the other hand, one can see easily that if (\ref{3-58}) is violated,  then $C\varrho$ remains an upper bound for $\rho$ at time $|t|\le T_0=Ch|\log h|\varsigma^{-1}$; therefore contribution of such elements to the remainder does not exceed
\begin{equation}
C\rho\gamma (\varsigma^2+\varrho\gamma+\mu h) h^{-4}T_0^* T_1^{-1} \asymp C\mu^{-1}h^{-3}\varrho\gamma  \varsigma^{-1}+
C\mu^{-1}h^{-3}\varrho^2\gamma^2 \varsigma^{-3}+
Ch^{-2}\rho\gamma  \varsigma^{-3}
\label{3-60}
\end{equation}
which does not exceed the same expression with $\varsigma=\gamma$ and corresponding $\varrho$; one can see easily that 
$\varrho\le \gamma h^{-5\delta'}$ and (\ref{3-60}) is $o(\mu^{-1}h^{-3})$.

Therefore only elements with $|\nabla 'V^{-1}f_j|\le C\gamma$  
remain to be treated, where either $|x_1|\le \gamma =C\mu^{-1}$ or 
$C\mu^{-1}\le |x_1|\asymp \gamma\le \mu^{-1+\delta}$.

\medskip
\noindent
(ii)
Let us consider the propagation speed with respect to $\rho$ more precisely. Note that as $|x_1|\asymp \gamma\ge C\mu^{-1}$ one can translate non-diagonal term $\mu^{-1}\Re (\omega Z_1^*Z_2^2)$ into 
\begin{equation*}
\mu^{-2}x_1^{-1}|\omega|^2 (Z_1^*Z_1-4z_2^*Z_2)Z_2^*Z_2+\dots
\end{equation*}
with 
\begin{equation*}
\rho= -\mu^{-2}x_2^{-1}|\omega|^2 (Z_1^*Z_1-4z_2^*Z_2)Z_2^*Z_2+\dots
\end{equation*}
Then 
$[a,\rho]=O(\mu^{-s}\gamma^{-s}+\mu^{-1})$ and furthermore along trajectories 
$[a,\rho](t)=[a,\rho](0)+O( (\mu^{-s}\gamma^{-s}|t|)$ (where $s$ is an arbitrarily large exponent and  $\delta=\delta(s)>0$ is small enough) and therefore 

\begin{claim}\label{3-61}
$\rho$ preserves both its sign and the bound from below as $\rho(0)\ge C\mu^{-s}\gamma^{1-s}$ and $|t|\le T_1$ with
\begin{equation}
T_1\Def \epsilon \rho^{1/2}\mu^{s/2}\gamma^{s/2}\qquad \text{as\ } {\bar\gamma}_0= C\mu^{-1}\le \gamma \le {\bar\gamma},\ \rho\ge C\mu^{-s/2}\gamma^{1-s/2}
\label{3-62}
\end{equation}
\end{claim}
(where $\rho=\rho (0)$). Therefore the contribution of the corresponding strip to the remainder does not exceed
\begin{equation*}
G(\rho,\gamma)\Def C\mu^2h^{-2}\gamma\times h\rho^{-2} \times \bigl((\rho^2+\gamma^2)(\mu h)^{-1}+1\bigr)\times \rho (\mu h)^{-1}\times \rho^{-1/2}\mu^{-s/2}\gamma^{-s/2}.
\end{equation*}
Thus after summation over $\rho\ge C\mu^{-s/6}\gamma^{1-s/6}$ \ I arrive to  $G(1,\gamma)+G(\mu^{-s/6}\gamma^{1-s/2})$ and for large enough $s$ summation with respect to $\gamma$ results in $G(1,\mu^{-1})+G(\mu^{-1}\mu^{-1})$ which is $O(\mu^{-1}h^{-3})$.

Meanwhile the contribution of the zone $\{|x_1|\asymp\gamma, \rho \le\varrho\Def  C\mu^{-s/6}\gamma^{1-s/6}\}$ to the remainder due to proposition \ref{prop-3-10} does not exceed $C\varrho \gamma \mu h^{-3}\asymp C\mu^{-s/6}\gamma^{1-s/6}h^{-3}$ and summation with respect to $\gamma$ results in $C\mu^{-1}h^{-3}$. Therefore
\begin{claim}\label{3-63}
Contribution of zone $\{{\bar\gamma}_0\le |x_1|\le {\bar\gamma}\}$ to the remainder is $O(\mu^{-1}h^{-3})$.
\end{claim}

\medskip
\noindent
(iii) Similar arguments work for zones 
$\{|x_1|\le {\bar\gamma}_0,\ \rho\ge \varrho\Def C\mu^{-1}\}$ with $T_1= C\rho ^{1/2}$
and for $\{|x_1|\le {\bar\gamma}_0,\ \rho\le \varrho\}$ and therefore 
Contribution of zone $\{ |x_1|\le {\bar\gamma}_0\}$ to the remainder is $O(\mu^{-1}h^{-3})$.

This concludes Part III and the whole proof.\end{proof}

\subsection{Sharp asymptotics at resonances}
\label{sect-3-5}

In this subsection I am going to prove sharp remainder estimate under generic assumptions to $V$. I know from proposition \ref{prop-1-6} that in the generic situation critical points of $\phi_\alpha=\alpha \log f_1+(1-\alpha)\log f_2 -\log V$ are non-degenerate except of discrete values of $\alpha=\alpha_j$. One can prove easily that

\begin{claim}\label{3-64}%\label{3N-51}%\label{3N-40}%\label{3C-65}
In the generic case degenerate critical points $\alpha_j$  of $\phi_\alpha$ are not resonances.  Then as $\alpha \ne \alpha_j$   the set of critical points is a smooth 1-dimensional curve parametrized by $\alpha$ and resonance surface is 3D surface.
\end{claim}
\begin{claim}\label{3-65}%\label{3N-52}
In the generic case these curve and resonance surface $\Xi_{kl}$ meet at isolated points and are transversal in them.
\end{claim}
Then
\begin{phantomequation}\label{3-66}\end{phantomequation}
\begin{multline}
{\tilde\nu} (\rho,\zeta,\gamma)\Def \mes \Bigl\{ (x,\alpha):\ |kf_1-lf_2|<\gamma,\\  |\nabla'\bigl( \alpha\log f_1+(1-\alpha)\log f_2- \log V\bigr)|\le \zeta,\\
\hskip83pt|\nabla \bigl( \alpha\log f_1+(1-\alpha)\log f_2- \log V\bigr)|\le \rho \Bigr\}\le C\rho^{q-r} \zeta^r\gamma  \\ 
\text{as\ }\ \rho \ge  \gamma,\ \rho\ge \zeta \ge (\rho\gamma)^{1/2}
\tag*{$(3.66)_{q,r}$}\label{3-66-q}
\end{multline}
with  $r=3$, $q=4$.

\begin{proposition}\label{prop-3-11} 
Let us assume that $|f_1-f_2|\ge\epsilon_0$ and  $|\nabla (f_1f_2^{-1})|\ge \epsilon_0$ on $\supp\psi$. Furthermore, let conditions $(\ref{0-8})$, \ref{3-21-q} with $q>3$ and \ref{3-66-q} be fulfilled with $r>2$,   $q>3$. Then  under condition $(\ref{3-18})$ asymptotics with the standard implicit main part $(\ref{0-9})$ and the remainder estimate $O(\mu^{-1}h^{-3})$ holds.
\end{proposition}

\begin{proof}[Proof, Part I]    I will follow the proof of proposition \ref{prop-3-7} and use the same numbering of its parts.  

\medskip
\noindent
(i) Repeating arguments leading to the proof of statements  of proposition \ref{prop-3-4} one can prove easily the following analogue of (\ref{3-40}):

\begin{claim}\label{3-67} Under condition \ref{3-21-q} the contribution of zone $\bigl\{x:\ |x_1| \ge \max(\rho, \mu^{-2/3})\bigr\}$ to the remainder does not exceed \ref{3-22-q}.
\end{claim}

\medskip
\noindent
(ii) Consider zone: $\bigl\{\rho \ge \max(\gamma,\mu^{-2}\gamma^{-2},\ \gamma\ge {\bar\gamma}\bigr\}$.
Then the total contribution of all $(\rho,\zeta,\gamma)$ elements does not exceed expression (\ref{3-38}) multiplied by ${\tilde\nu}(\rho,\zeta,\gamma)\rho^{-1}$:
\begin{equation}
C\mu h^{-1}\bigl(\zeta (\mu h)^{-1}+\zeta^{-1}\bigr)\times \bigl(\rho (\mu h)^{-1}+1\bigr) \times \rho^{q-r-1}\zeta^r
\label{3-68}
\end{equation}
where I used \ref{3-66-q} to estimate ${\tilde\nu}$. Then as $r>1$ summation with respect to $\zeta$ results in the same expression as $\zeta= \rho$;
further as $q>3$ summation with respect to $\rho$ results in the same expression as $\rho=1$ which is $C\mu^{-1}h^{-3}$ and summation with respect to $\gamma$ results in $C\mu^{-1}h^{-3}|\log h|$. One can get rid of the logarithmic factor using the same arguments as in the Part I (ii) of proof of proposition \ref{prop-3-4}.

\medskip
\noindent
(iii) Consider the remaining part of zone $\{\epsilon \ge\gamma\ge{\bar\gamma}\}$ and introduce scaling function $\eta$ by (\ref{3-41}). Then one gets \ref{3-37-*} modified in the same way as in the Part I (iii) of proof of proposition \ref{prop-3-4} and  multiplied by $\eta\gamma^{-1} {\tilde\nu}(\rho\gamma)^{-1}$:
\begin{equation}
C\mu h^{-1} 
\bigl(\zeta  (\mu h)^{-1}+\zeta^{-1}\bigr)\times \bigl(\rho (\mu h)^{-1}+1\bigr)
\times \eta\gamma^{-1}\times \rho^{q-r-s-1}\zeta^r.
\tag*{$(3.64)^*$}
\label{3-68-*}
\end{equation}
Then as in (ii) summation with respect to $\zeta$ results in its value as $\zeta=\rho$:
\begin{equation*}
C\mu h^{-1} 
\bigl(\rho  (\mu h)^{-1}+\rho^{-1}\bigr)\times \bigl(\rho (\mu h)^{-1}+1\bigr)
\times \mu^2\varrho^{-2}\gamma^3\rho\gamma^{-1}\times \rho^{q-s-1}\gamma^{s-1}
\end{equation*}
and summation with respect to $\rho$ returns the above expression at its largest value which is $\mu^{-2}\gamma^{-2}\varrho^2$; the result does not exceed 
$C\mu^{-1-\delta'}h^{-3}$ and summation with respect to $\gamma$ returns $o(\mu^{-1}h^{-3})$,
\end{proof}

\begin{remark}\label{rem-3-12} t
(i) Again as the order of resonance $m\ge 4$, analysis of (iii) is not needed;

\medskip
\noindent
(iii) Furthermore,  the rough remainder estimate of zone $\{|x_1|\le {\bar\gamma}\}$ returns $O(\mu^{-1}h^{-3})$ as $m\ge5$ and $\mu\le h^{\delta'-1}$ and $O(\mu^{-1+\delta}h^{-3})$ as \underline{either} $m=4$  
\underline{or} $m=5$, $h^{\delta'-1}\le \mu\le ch^{-1}$.
\end{remark}

\begin{proof}[Proof, Part II]\label{pf-3-11-II} (i) Analysis in zone $\{|x_1|\le {\bar\gamma}\}$ is now  simpler.   Note first that the contribution of zone $\{\zeta\le {\bar\gamma}\}$ to the remainder does not exceed 
$C\mu h^{-3}{\bar\gamma}^{r+1}=O(\mu^{-1}h^{-3})$  as $r>1$ and $\delta<\delta(r)$.

\medskip
\noindent 
(ii) Consider zone $\{\zeta \ge C{\bar\gamma}\}$. Then defining $T_0,T_0^*$ by (\ref{3-35}) and $T_1$ by (\ref{3-36}) one estimates contribution of $(\rho,{\bar\gamma},\zeta)$ elements by (\ref{3-52}) multiplied by ${\tilde\nu}(\rho{\bar\gamma})^{-1}$:
\begin{equation}
C\mu^2 h^{-2}{\bar\gamma} \times h\zeta^{-2}\times 
\bigl(\zeta^2 (\mu h)^{-1}+ 1\bigr)\times \rho (\mu h)^{-1} \times 
\bigl( \mu^{-1} {\bar\gamma} ^{-1/2}+\rho^{-1}\bigr)\times \rho^{q-r-2}\zeta^r
\label{3-69}
\end{equation}
and summation with respect to $\zeta$ results in its value as $\zeta=\rho$ (now I need $r>2$)
\begin{equation*}
C\mu  h^{-2}{\bar\gamma} \times  
\bigl(\rho^2 (\mu h)^{-1}+ 1\bigr)\times    
\bigl( \mu^{-1} {\bar\gamma} ^{-1/2}+\rho^{-1}\bigr)\times \rho^{q-3}
\end{equation*}
and summation with respect to $\rho$ results in $Ch^{-3}{\bar\gamma}$ which is marginally worse than $O(\mu^{-1}h^{-3})$. To improve it one can sum to $\zeta \le \mu^{-\kappa}$ and in zone $\{\zeta\ge \mu^{-\kappa}\}$ one can take $T_1=\epsilon \mu\zeta$.
\end{proof}

\subsection{General estimates near $\Sigma$}
\label{sect-3-6}

Now let us consider the vicinity of $\Sigma=\{f_1=f_2\}=\{v_1=v_2=0\}$ where
\begin{claim}\label{3-70}
$v_1|_\Sigma=v_2|_\Sigma=0$,  $(\nabla v_1)|_\Sigma$ and $(\nabla v_2)|_\Sigma$ are linearly independent
 \end{claim} 
 and near $\Sigma$
\begin{equation}
f_{1,2}=f\pm (v_1^2+v_2^2)^{1/2},\qquad f>0.
\label{3-71}
\end{equation}
This analysis is simpler than near third  order resonances because codimension of $\Sigma$ is 2 and everywhere factor $\gamma^1$ reflecting measure should be replaced by $\gamma^2$. 

\begin{proposition}\label{prop-3-13} Let condition $(\ref{1-1})-(\ref{1-2})$ be fulfilled and let $\psi$ be  supported in the small vicinity of $\Sigma$.

The standard implicit asymptotic formula $(\ref{0-9})$ holds with the remainder estimate $O(\mu^{-1}h^{-3}+\mu^2h^{-2})$.
\end{proposition}

\begin{proof} (i) Note first that

\begin{claim}\label{3-72}%\label{3N-55}
The contribution of zone $\{\dist (x,\Sigma)\le \gamma\}$ to the remainder does not exceed $C\mu h^{-3}\gamma^2$ and as 
\begin{equation}
\gamma\Def\dist (x,\Sigma)\asymp  |f_1-f_2|\le {\bar\gamma}_1\Def c\mu^{-1}+ c(\mu h)^{1/2}
\label{3-73}%\label{3N-56}
\end{equation}
this contribution is $O(\mu^{-1}h^{-3}+\mu^2h^{-2})$.
\end{claim}

On the other hand, 

\begin{claim}\label{3-74}
One can reduce operator to the canonical form without non-diagonal terms as long as
\begin{equation}
\gamma\ge {\bar\gamma}\Def  \mu^{-1/2}h^{1/2-\delta'}+\mu^{-2}h^{-\delta'}.
\label{3-75}
\end{equation}
\end{claim}

The second term in (\ref{3-75})  appears because one  needs to get rid off terms $Z_1^iZ_1^{*\,j}Z_2^kZ_2^{*\,l}$ with $i+k=j+l$ but $(i,j)\ne (k,l)$ and these terms are of magnitude $O(\mu^{-2})$ unless $i+j+k+l=2$ in which case one just diagonalizes the quadratic form and it is where the first term in in (\ref{3-75})  comes from.

Important is that ${\bar\gamma}\le {\bar\gamma}_1$. In the quest for remainder estimate $O(\mu^{-1}h^{-3})$ one would need to take ${\bar\gamma}= c\mu^{-1}$ and ${\bar\gamma}\le {\bar\gamma}_1$ would hold as $\mu \le h^{\delta'-1}$.

\medskip
\noindent 
(ii) Making $\epsilon\gamma$-admissible partition and $\rho$-admissible subpartition with 
\begin{equation}
\rho =\epsilon |\nabla \phi_\alpha|\gamma + {\frac 1 2}{\bar\rho} ,\qquad {\bar\rho} =(C\mu^{-1}h|\log h| )^{1/2}
\label{3-76}
\end{equation}
one can take 
\begin{equation}
T_0^*= Ch\rho^{-2}\gamma,\qquad T_1=\epsilon \mu\gamma
\label{3-77}
\end{equation}
and  the total contribution to the remainder of $(\gamma,\rho)$ subelements with $\rho\ge \varrho $ to the remainder does not exceed
\begin{multline}
C\mu^2h^{-2}\times \mu^{-1}h\rho^{-2} \times \bigl(\rho^2(\mu h \gamma)^{-1}+1\bigr)\times \bigl(\rho (\mu h\gamma)^{-1}+1\bigr)\times \gamma^2\asymp\\
C\mu^{-1}h^{-3}\rho  + C\rho^{-1}h^{-2}\gamma  + C \mu h^{-1} \rho^{-2}\gamma^2
\label{3-78}
\end{multline}
where $\gamma^2$ is their total measure. The right-hand expression sums with respect to $\rho\in({\bar\rho} ,\gamma) $ to 
$G(\gamma)\Def C\mu^{-1}h^{-3}\gamma  + Ch^{-2}{\bar\rho} ^{-1} \gamma+ C\mu h^{-1}{\bar\rho}^{-2}\gamma^2 $. Then with respect to $\gamma$   it sums to  
$G(1)=O\bigl( \mu^{-1}h^{-3}+ \mu^{1/2}h^{-5/2} +\mu^2h^{-2}\bigr)$ where the middle term is less than the sum of two others.

\medskip
\noindent 
(iii) Meanwhile the total contribution to the remainder of $({\bar\rho} ,\gamma)$ subelements does not exceed
\begin{multline}
C\mu^2h^{-2}\times   \bigl({\bar\rho} ^2(\mu h \gamma)^{-1}+1\bigr)\times \bigl({\bar\rho}  (\mu h\gamma)^{-1}+1\bigr)\times \gamma^2\asymp\\
Ch^{-4}{\bar\rho}^3  +C\mu h^{-3}\gamma {\bar\rho}  + C\mu^2h^{-2}\gamma^2.
\label{3-79}
\end{multline}
This expression sums with respect to $\gamma$ to its value as $\gamma=1$ resulting in  
\begin{equation*}
C\mu^{-3/2} h^{-5/2}|\log h|^{3/2}+C\mu^{1/2}h^{-5/2}|\log h|^{1/2}+C\mu^2h^{-2}.
\end{equation*}
Note that the first and the third term are properly estimated and the second terms is properly estimated  save border case $h^{-1/3}|\log h|^{-K}\le \mu \le h^{-1/3}|\log h|^K$ which is treated as in the part (iii) of the proof of proposition \ref{prop-3-1}.

\medskip
\noindent 
(iv)
Finally, as $\gamma\le {\bar\gamma}$ one does not need a subpartition; the contribution to the remainder does not exceed $C\mu h^{-3}{\bar\gamma}^2$ due to an analogue  of proposition \ref{prop-3-10} below: 
\end{proof}

\begin{proposition}\label{prop-3-14} Proposition \ref{prop-3-10} remains true near $\Sigma$ (i.e. without condition $|f_1-f_2|\ge \epsilon$ provided at point ${\bar x}$ main part of precanonical form is 
$f_1(\mu^2x_3^2+h^2D_3^2)+f_2(\mu^2x_4^2+h^2D_4^2)$.
\end{proposition}

\begin{proof}
Proof basically repeats the one of proposition \ref{prop-3-10}.
\end{proof}

\subsection{Sharp asymptotics near $\Sigma$}
\label{sect-3-7}

\subsubsection{}
\label{sect-3-7-1}
Now let us improve the results of the previous subsection. Let us note that
\begin{claim}\label{3-80}%\label{3N-67}
Microhyperbolicity condition holds at ${\bar x}\in \Sigma$ iff in frames of 
$\nabla (fV^{-1})$ is not a linear combination of $\nabla (v_1V^{-1})$, 
$\nabla (v_2V^{-1})$ with coefficients $(\beta_1,\beta_2)\in \bR^2\cap B(0,1)$ (\footnote{\label{foot-14} In the uniform sense, i.e. (\ref{1-22})}.
\end{claim}

In virtue of \cite{IRO4} I have already

\begin{proposition}\label{prop-3-15}
If $(\ref{1-22})$ is fulfilled on $\supp \psi$   then the standard formula holds with ${\bar T}=Ch|\log h|$.
\end{proposition}

Now let us analyze the meaning of (\ref{1-22}). First of all, it is fulfilled as 
$\nabla_\Sigma (fV^{-1})\ne 0$. So one needs to consider only set $\Sigma_0$ of the critical points of $fV^{-1}\bigr|_\Sigma$:
\begin{equation}
\Sigma_0=\bigl\{ x\in \Sigma,\ \nabla_\Sigma (fV^{-1})=0\bigr\}.
\label{3-81}
\end{equation}
Then 
\begin{proposition}\label{prop-3-16}
For generic $V$ 
\begin{claim}\label{3-82}%\label{3N-69}
 $\Sigma_0$ consists of separate non-degenerate points;
 \end{claim}
\begin{claim}\label{3-83}%+5
Magnetic form $\omega_F$ restricted to $\Sigma$ is the generic closed form on $\Sigma$ and thus degenerates on the smooth curve $\{ \{v_1,v_2\}=0\}$;
\end{claim}
\begin{claim}\label{3-84}%+6
$\omega_M $ does not degenerate on $\Sigma_0$ (which is equivalent to $\{v_1,v_2\}\ne 0$ on $\Sigma_0$).
\end{claim}
\end{proposition}

One can write down many generic properties, but they  are overkill.

\subsubsection{}
\label{sect-3-7-2} 
First let us improve the remainder under condition (\ref{3-85}) below (I remind that $f={\frac 1 2}(f_1+f_2)$):

\begin{proposition}\label{prop-3-17}
Let at some point ${\bar x}\in \Sigma$
\begin{equation}
|\nabla (fV^{-1})|\ge \epsilon_0.
\label{3-85}
\end{equation}
Then as $\psi$ is supported in the small enough vicinity of ${\bar x}$
the standard implicit formula holds with  remainder  
$O\bigl(\mu^{-1}h^{-3}+ \mu^{3/2}h^{-3/2}|\log h|\bigr)$.
\end{proposition}

\begin{proof} (i) Note, that here in contrast to (\ref{1-22}) the complete gradient is considered. One needs to consider the case when microhyperbolicity condition is not fulfilled in ${\bar x}=0$; then
\begin{equation}
\nabla  (fV^{-1}) = \beta_1 \nabla (v_1V^{-1}) +\beta_2\nabla (v_2V^{-1})\qquad \text{at\ } {\bar x}, \quad  \beta =(\beta_1,\beta_2)\in \bR^2: 1\ge |\beta|\ge \epsilon_0
\label{3-86}
\end{equation} 
where $|\beta|\ge\epsilon_0$ due to (\ref{3-85}).

Let us consider $\phi_\alpha=\alpha (f_1V^{-1})+(1-\alpha) (f_2V^{-1})$; note that one can extend $\beta_j$ to the vicinity of ${\bar x}$ so that  \begin{equation*}
fV^{-1}={\tilde f}(w)+\beta_1 (w) v_1V^{-1}+\beta_2(w) v_2V^{-1}+O(|v|^2)
\end{equation*}
where $w=(w_1,w_2)$ are coordinates on $\Sigma$.  

Without any loss of the generality one can assume that $\beta_2(w)=0$; then 
under  conditions (\ref{3-85}), (\ref{3-86}) 
\begin{equation}
|\nabla \phi_\alpha|\ge \epsilon_1
\bigl(|(v_2-\omega v_1^2)(v_1^2+v_2^2)^{-1/2}|+ |\alpha-{\bar\alpha}|\bigr)
\label{3-87}
\end{equation}
where $\omega=\omega(w,v_1)$ is a smooth function.

\medskip
\noindent
(ii) Let us follow  the proof of proposition \ref{prop-3-13}. In part (ii) estimate (\ref{3-78}) (for contribution of all elements with $\rho\ge\varrho$) gains a  factor $\rho\gamma^{-1}$ and becomes
\begin{equation}
C\mu^{-1}h^{-3}\rho^2\gamma^{-1}  + Ch^{-2} + C \mu h^{-1}\rho^{-1}\gamma.
\tag*{$(3.78)^*$}\label{3-78-*}
\end{equation}
This expression sums with respect to $\rho$ ranging from ${\bar\rho}$ to $\gamma$ to
\begin{equation*}
C\mu^{-1}h^{-3} \gamma  + C  h^{-2} |\log h| + C \mu h^{-1} {\bar\rho}^{-1}\gamma.
\end{equation*}
Then summation with respect to $\gamma$ results in 
\begin{equation*}
C\mu^{-1}h^{-3}  + C h^{-2} |\log h|^2 + C \mu ^{3/2}h^{-3/2} 
\end{equation*}
which is less than the announced estimate.

\medskip
\noindent 
(iii) In part (ii) estimate (\ref{3-79}) (for contribution of all elements with $\rho\le{\bar\rho}$) gains a  factor ${\bar\rho}\gamma^{-1}$ and becomes
\begin{equation}
Ch^{-4}{\bar\rho}^4\gamma^{-1}  +C\mu h^{-3} {\bar\rho}^2  + C\mu^2h^{-2}{\bar\rho}\gamma.
\tag*{$(3.79)^*$}\label{3-79-*}
\end{equation}
This expression sums with respect to $\gamma$ to 
\begin{equation*}
C\mu^{-2}h^{-2}|\log h|^2{\bar\gamma}^{-1}+C h^{-2}|\log h|^2+ C\mu^{3/2}h^{-3/2}|\log h|^{1/2}
\end{equation*}
which is also below than announced remainder estimate.

\medskip
\noindent 
(iv)
Finally, as $\gamma\le {\bar\gamma}$ one does not need a subpartition; the contribution to the remainder does not exceed $C\mu h^{-3}{\bar\gamma}^2= C\mu^{-1}h^{-3} + C\mu^{-1/2}h^{\delta-5/2}$ and the second term here is obviously much less than the announced estimate.
\end{proof}

So, condition (\ref{3-85}) is a kind of non-degeneracy condition, improving remainder estimate.

\subsubsection{}
\label{sect-3-7-3} 
Let us now use the remaining non-degeneracy conditions.

\begin{proposition}\label{prop-3-18} Let conditions $(\ref{3-18})$ and
\begin{phantomequation}\label{3-88}\end{phantomequation}
\begin{multline}
|(2p+2)\mu h f V^{-1} -1|+
|\nabla_\Sigma\bigl( (2p+2)\mu h f V^{-1}-1\bigr)|\le \epsilon_0 \implies\\ \Hess _\Sigma \bigl((2p+2)\mu h f V^{-1}-1\bigr)  \text{\ \ has at least $r$ eigenvalues  }\\ \text{with absolute values greater than $\epsilon_0$}
\tag*{$(3.88)_r$}\label{3-88-r}
\end{multline}
be fulfilled. Then

\medskip
\noindent 
(i) The total  remainder is given by $(\ref{3-21})_{r+1}$ while the main part is given by the standard implicit formula $(\ref{0-9})$. 

\medskip
\noindent 
(ii) Under condition $(\ref{3-85})$ the total  remainder is given by $(\ref{3-21})_{r+2}$ as $r=1$ while the main part is given by the standard implicit formula $(\ref{0-9})$;

\medskip
\noindent 
(iii) Under conditions $(\ref{3-85})$ and $(\ref{3-84})$ the total  remainder is $O(\mu^{-1}h^{-3})$ as $r=2$ while the main part is given by the standard implicit formula $(\ref{0-9})$.
\end{proposition}

\begin{proof} (a) Under condition (\ref{3-85}) the total contribution to the remainder of all $(\gamma,\rho)$-elements does not exceed \ref{3-78-*} i.e.
\begin{phantomequation}\label{3-89}\end{phantomequation}
\begin{equation}
C \mu^{-1}h^{-3}\rho ^{k+1}\gamma^{-k} + C\rho^{k-1}\gamma^{1-k} h^{-2} +  C\mu h^{-1}\gamma^{2-k} \rho^{k-2}
\tag*{$(3.89)_k$}\label{3-89-k}
\end{equation}
with $k=1$ while without it it does not exceed the same expression with $k=0$.

Further, under extra condition \ref{3-88-r} this expression acquires  factor $C\rho^r\gamma^{-r}$. Therefore again $(\ref{3-89})_q$  with $q=k+r$ gives a proper estimate for the total contribution of all $(\gamma,\rho)$-elements to the remainder. 

Consider summation with respect to ${\bar\rho}\le \rho \le\gamma\le 1$. The first term sums to $C\mu^{-1}h^{-3}$ independently on $q$; the second term sums to $C h^{-2}|\log h|^2$ as $q=1$ and to $Ch^{-2}|\log h|$ as $q\ge 2$; the third term sums to $C\mu h^{-1}{\bar\rho}^{-1}$ as $q=1$, to $C\mu h^{-1}|\log h|^2$ as $q=2$ and to $C\mu h^{-1}|\log h|$ as $q=3$. 

Therefore in all cases but one the remainder estimate \ref{3-21-q} with $q=k+r$ is proven; this exceptional case is $q=3$, $\mu \ge h^{-1}|\log h|^{-1}$ when the  estimate $O(h^{-2}|\log h|)$ is recovered; I remind that only zone $\{\rho\ge {\bar\rho},\ \gamma\ge {\bar\gamma}\}$ is covered so far.

\medskip
\noindent
(b) To cover the remaining case let us introduce a scaling function $\zeta=\epsilon |\nabla_\Sigma fV^{-1}|+{\bar\rho}$ on $\Sigma$; one can extend this function to
\begin{equation}
\zeta=\epsilon |\nabla_\Sigma fV^{-1}|+\gamma.
\label{3-90}
\end{equation}
Then in arguments above one can replace factor $\rho^2\gamma^{-2}$ by $\zeta^2$ and then contribution of zone $\{\zeta\le\gamma^{\delta'}\}$ would be $O(\mu^{-1}h^{-3})$ as $q=3$; so only subzone $\{\zeta\ge\gamma^{\delta'}\}$ remains to be treated. In this subzone one can take
\begin{equation}
T_0^*= Ch \zeta^{-2},\qquad T_1= \epsilon \mu \gamma
\label{3-91}
\end{equation}
where one can take $T_0^*= Ch \zeta^{-2}$ rather than $T_0^*= Ch \zeta^{-1}\gamma^{-1}$ due to (\ref{3-84}). Then the total contribution of such $(\gamma,\zeta)$ elements to the remainder does not exceed
\begin{equation*}
C\mu^2h^{-2}\zeta^2\gamma^2 \times h\zeta^{-2}\times \mu^{-1}\zeta^{-1}\times \bigl( (\zeta^2+\gamma)(\mu h)^{-1} +1\bigr) \times (\mu h)^{-1}
\end{equation*}
where $\zeta^2\gamma^2$ is the measure and $(\mu h)^{-1}$ and 
$\bigl( (\zeta^2+\gamma)(\mu h)^{-1} +1\bigr) $ are  estimates of the numbers of indices ``$p$'' and corresponding indices ``$n$''; noting that $\zeta^2\ge \gamma$ one can rewrite this expression as
\begin{equation*}
C  \mu^{-1}h^{-3}\gamma^2   \zeta + C  h^{-2}\gamma^2    \zeta^{-1}
\end{equation*}
and the summation with respect to $\zeta\ge \gamma^{\delta'}$ and $\gamma$ results in  $O(  \mu^{-1}h^{-3})$.

Thus estimate \ref{3-21-q} for contribution of zone $\{\rho\ge {\bar\rho},\ \gamma\ge {\bar\gamma}\}$ is established.

\medskip
\noindent
(c) Further, contribution of  $(\gamma,\rho)$ elements with $\rho\le {\bar\rho}$, $\gamma\ge {\bar\gamma}$ to the remainder does not exceed (\ref{3-78}) as $k=0$ or \ref{3-78-*} as $k=1$ multiplied by ${\bar\rho}^r\gamma^{1-r}$ i.e.
\begin{phantomequation}\label{3-92}\end{phantomequation}
\begin{equation}
Ch^{-4}{\bar\rho} ^{q+3}\gamma^{-q} + C\mu h^{-3} {\bar\rho}^{q+1}\gamma^{1-q} + C \mu ^2h^{-2} {\bar\rho}^q\gamma^{3-q}
\tag*{$(3.92)_q$}\label{3-92-q}
\end{equation}
and one can check easily that summation with respect to $\gamma\ge {\bar\gamma}$ results in expression not exceeding $(\ref{3-19})_q$. 

\medskip
\noindent
(d) Finally, contribution of zone $\{\gamma\le {\bar\gamma}\}$ to the remainder does not exceed $C\mu h^{-3}{\bar\gamma}^2$ which is $O(\mu^{-1}h^{-3})$ unless $\mu \ge h^{2\delta-1}$ and even in this case it is less than $(\ref{3-19})_1$. 

Therefore as $q\ge 2$ I need some better arguments in this zone. Again, one needs to consider only part of it with $\{\zeta\ge \mu^{-\delta'}\}$ as contribution of zone $\{\zeta\le {\bar\zeta}\}$ does not exceed $C\mu h^{-3}{\bar\gamma}^2{\bar\zeta}^2$.

As using precanonical form the speed would be $O(\mu^{-1})$ and since only a magnitude of $\zeta$ is important for us and also inequality $\zeta\ge \gamma^{\delta'}$, one can take
\begin{equation}
T_0^*= Ch \zeta^{-1}{\bar\gamma}^{-1},\qquad T_1= \epsilon  \mu {\bar\gamma} ^{1-\delta''} 
\label{3-93}
\end{equation}
where an appropriate time direction for this $T_1$ is taken. In these arguments I do not assume (\ref{3-84}) and and thus $T_0$ is not as it was in (b) (surely, some improvements are possible but not needed). 

Then the total contribution to the remainder of $\zeta$-elements does not exceed
\begin{equation*}
Ch^{-4}\zeta^r {\bar\gamma}^2 \times h\zeta^{-1}{\bar\gamma}^{-1}\times \mu^{-1}{\bar\gamma}^{\delta''-1} \asymp C\mu ^{-1}h^{-3}\zeta^{r-1}{\bar\gamma}^{\delta''}
\end{equation*}
which sums with respect to $\zeta$ to  $C\mu^{-1}h^{-3}|\log h|{\bar\gamma}^{\delta''}=o(\mu^{-1}h^{-3})$ even as $r=1$.
\end{proof}

\begin{remark}\label{rem-3-19} (i) Since in the generic case in the critical points of $fV^{-1}|_\Sigma$ are non-degenerate (i.e. (\ref{3-82})   is fulfilled) and also (\ref{3-83}),(\ref{3-84}) hold, the remainder estimate is $O(\mu^{-1}h^{-3})$.

\medskip
\noindent 
(ii) I think one can get rid off logarithmic factors in the above estimates but I do not care.
\end{remark}

\subsection{Summary}
\label{sect-3-8}

\begin{proposition}\label{prop-3-20} {\rm (I)} Let $(g^{jk})$ be fixed and then $(V_j)$ be   be generic, more precisely:

\medskip
\noindent
{\rm (i)} Outside of $\Sigma=\{x:\ f_1=f_2\}$ critical points of $f_1f_2^{-1}$ satisfy $(\ref{2-30})_3$ and $(\ref{2-2})$;

\medskip
\noindent
{\rm (ii)} $\Sigma$ be smooth 2-dimensional manifold and $|f_1-f_2|\asymp \dist(x,\Sigma)$.

\medskip
\noindent
 Let  $V$ be general but satisfying $(\ref{0-8})$ at $\supp\psi$. Then under condition $(\ref{3-1})$  the remainder is $O(\mu^{-1}h^{-3}+\mu^2h^{-2})$ while the main part of asymptotics is given by  implicit formula $(\ref{0-9})$.

\medskip
\noindent
{\rm (II)} Furthermore, let $V$ be generic i.e. 

\medskip
\noindent
{\rm (iii)} Outside of $\Sigma$ and resonances condition $(\ref{3-21}_4$ be fulfilled;

\medskip
\noindent
{\rm (iv)} Near resonances condition $(\ref{3-66})_{4,3}$ be fulfilled;

\medskip
\noindent
{\rm (v)} At $\Sigma$ conditions $(\ref{3-82})-(\ref{3-85})$ be fulfilled\footnote{\label{foot-15} Actually these conditions should be fulfilled at $\Sigma_0^+\cup\Sigma_0^0$ only.}.

\medskip
\noindent Let  $V$ satisfy $(\ref{0-8})$ at $\supp\psi$. Then under condition $(\ref{3-18})$  the remainder is $O(\mu^{-1}h^{-3})$ while the main part of asymptotics is given by   implicit formula $(\ref{0-9})$.
\end{proposition}

\section{Calculations}
\label{sect-4}

The purpose of this section is to pass from the implicit formula (\ref{0-9}) to more explicit one, namely either (\ref{0-7}) or (\ref{0-5}) with more or less explicit expression for $\cE^\MW_\corr$. It will be done by different methods depending on the magnitude of $\mu$ and also non-degeneracy conditions and the methods applied will be used in the classification. My main concern will be either no non-degeneracy condition or the generic non-degeneracy condition.

\subsection{Temperate magnetic field}
\label{sect-4-1}

In this subsection I will use  formula (\ref{0-9}) with $T_0\le \epsilon \mu^{-1}$ and derive a remainder estimate. More precisely, if in section \ref{sect-3} $T$ was given by (\ref{0-9}) with $T_0\ge \epsilon \mu ^{-1}$ I replace it by $T_0=\epsilon \mu^{-1}$ and estimate an error. 

This error actually is the contribution of the affected domain to the remainder with some $T_0\ll \mu^{-1}$ (usually $T_0=Ch|\log h|$ under condition (\ref{0-8})) and $T_1=\epsilon \mu^{-1}$. 

Further, under condition (\ref{0-8}) I can trade $T \le \epsilon \mu^{-1}$ to $T= Ch|\log h|$ and then I can apply the standard approach; I will assume here by default that 
\begin{equation}
h^{-\delta_0}\le \mu \le h^{-1+\delta_0}
\label{4-1}
\end{equation}
leaving case $\le h^{-1+\delta_0}\le \mu \le ch^{-1}$ for a separate consideration.

\subsubsection{}
\label{sect-4-1-1} Let us consider regular points first. 
 
\begin{proposition}\label{prop-4-1}  Assume that  there are no resonances of order not exceeding $M$ and also $|\nabla (f_1f_2^{-1})|\ge \epsilon_0$  on $\supp\psi$. Further, assume that condition \ref{3-21-q} with $q\ge1$ is fulfilled. Then  under condition $(\ref{4-1})$ the remainder estimate is given by 
\begin{phantomequation}\label{4-2}\end{phantomequation}
\begin{equation}
O\Bigl(\mu^{-1}h^{-3}+ (\mu h)^{(q+2)/2}h^{-4}|\log h|^{q/2}\Bigr)
\tag*{$(4.2)_q$}
\label{4-2-q}
\end{equation}
while the main term of asymptotics given by $(\ref{0-9})$ with any $T\ge \epsilon \mu^{-1}$.
\end{proposition}

\begin{proof} First of all, as $q=1$ this remainder estimate is $O\bigl(\mu^{-1}h^{-3}+ (\mu h)^{3/2}h^{-4}|\log h|^{1/2}\bigr)$ which is no smaller than $O(\mu^{-1}h^{-3}+\mu^2h^{-2})$ given by proposition \ref{prop-3-1} and as $q>1$ remainder estimate \ref{4-2-q}  is no smaller than  \ref{3-22-q} given by proposition \ref{prop-3-4}; so both of these propositions could be applied and one needs to estimate a substitution $T_0\mapsto \epsilon \mu^{-1}$ error.

According to the proofs of these propositions $T_0=Ch \rho^{-2}|\log h|$ which is less than $\epsilon \mu^{-1}$  as 
\begin{equation}
\rho \ge {\bar\rho}_1\Def (C\mu h|\log h|)^{1/2};
\label{4-3}
\end{equation}
So one can take there $T_0=\epsilon \mu^{-1}$.

On the other hand, contribution of zone $\{\rho \le {\bar\rho}_1\}$ to the remainder with $T=\epsilon \mu^{-1}$ does no exceed $C\mu h^{-3}\nu ({\bar\rho}_1)\le C\mu h^{-3}{\bar\rho}_1^q$ which is exactly the second term in \ref{4-2-q}.
\end{proof}

\begin{corollary}\label{cor-4-2}  Let conditions of  proposition \ref{prop-4-1} be fulfilled. Then under assumption $(\ref{0-8})$  

\medskip
\noindent
{\rm (i)} The remainder estimate is given by \ref{4-2-q} while the main term of asymptotics given by 
\begin{equation}
\int \cE^\MW(x,0)\psi (x)\, dx.
\label{4-4}
\end{equation}

\medskip
\noindent
{\rm (ii)} In particular  the remainder estimate is $O(\mu^{-1}h^{-3})$ as $\mu \le {\bar\mu}_q\Def C(h|\log h|)^{-q/(q+4)}$ where in the general case 
${\bar\mu}_1= C(h|\log h|)^{-1/5}$ and in the generic case ${\bar\mu}_4= C(h|\log h|)^{-1/2}$.
\end{corollary}

\begin{proof} According to proposition \ref{prop-4-1} with the remainder estimate \ref{4-2-q} the main part of asymptotics is given by (\ref{0-9}) with $T=\epsilon \mu^{-1}$. However then under condition (\ref{0-8}) one can trade $T=\epsilon \mu^{-1}$ to any $T\ge {\bar T}\Def Ch|\log h|$. The rest is proven by the standard successive approximation method, applied to the original operator (rather than to the canonical form); one can take as an unperturbed operator the same operator with the coefficients frozen at $y$ which leads to the Weyl expression perturbed by $\sum_{m\ge 0, n\ge 1}\kappa_{mn}h^{-4+2m+2n}\mu^{2n}$ where terms with $n\ge 2$ or $m\ge 1$ do not exceed the announced remainder estimate. Alternatively one can take as an unperturbed operator the same operator with $g^{jk}$, $V$ frozen at $y$ and with $V_j$ replaced by $V_j(y)+\sum_k (\partial_kV_j)(y)(x_k-y_k)$ which leads directly to Magnetic Weyl expression. Details see in (\ref{IRO3-1-11}) of \cite{IRO3}.
\end{proof}

\subsubsection{}
\label{sect-4-1-2}
Now I want to improve the result in the general case and allow non-degenerate critical points of $f_1f_2^{-1}$.

\begin{proposition}\label{prop-4-3}  Assume that $f_1\ne f_2$ and there are no resonances of order not exceeding $M$  and also critical points of $f_1f_2^{-1}$  are non-degenerate on $\supp\psi$. Then under condition $(\ref{4-1})$ the remainder estimate is given by \ref{4-2-q} with  any $q<2$ arbitrarily close to $2$ while the main term of asymptotics given by the standard implicit formula with any $T\ge \epsilon \mu^{-1}$.
\end{proposition}

\begin{proof} Again, this remainder estimate is no smaller than $O(\mu^{-1}h^{-3}+\mu^2h^{-2})$ given by proposition \ref{prop-3-2}, so again one needs to estimate an  error arising from the substitution $T_0\mapsto \epsilon \mu^{-1}$.

Proof follows ideas of the proof of proposition \ref{IRO8-prop-4-5} of \cite{IRO8}. Let us introduce functions $\ell_k$, $k=1,\dots,K$  in the same proof but with $\gamma=1$ and again let ${\bar\ell}_k= (\mu h|\log h|^J)^{1/(k+1)}$. 

\medskip
\noindent
(i) Assume first that $|\nabla (f_1f_2^{-1})|\ge \epsilon_0$. Then  the arguments of the mentioned proof survive with this simplification. 

\medskip
\noindent
(ii) Assume now that $f_1f_2^{-1}$ has one non-degenerate critical point ${\bar x}$.   Consider first zone where $\ell_K(x)\le \gamma(x) \Def {\frac 1 2}|x-{\bar x}|$. 

Then the number of indices ``$p$'' does not exceed 
$C\bigl(\ell_k^K(\mu h \gamma)^{-1}+1\bigr)$ (where the second term is the smallest one) while the number of indices ``$n$'' for each $p$ does not exceed $C\bigl( \ell_K^{K+1}(\mu h)^{-1}+1\bigr)$ and the total contribution to the asymptotics of all such elements with $\ell_K\asymp {\bar\ell}_K$ and fixed magnitude of $\gamma(x)\asymp \gamma$ for some $k$  does not exceed 
\begin{equation*}
C\mu^2h^{-2}h^{-2}\times \bigl( \ell_K^{K+1}(\mu h)^{-1}+1\bigr) \times \ell_k^K(\mu h \gamma)^{-1} \times \gamma^4 
\end{equation*}
where  $\gamma^4$ is the total measure of zone $\{\gamma(x)\asymp \gamma\}$.
This expression does not exceed $C\mu^2h^{-2}{\bar\ell}_K^{-1}\gamma^3|\log h|^J$ and summation with respect to $\gamma$ results in 
$C\mu^2h^{-2}{\bar\ell}_K^{-1}|\log h|^J\lesssim C(\mu h)^{2-\delta}h^{-4}$ as
$\delta>1/K$.

On the other hand, repeating arguments of the proof of proposition \ref{IRO8-prop-4-5} of \cite{IRO8} one can see easily that the total contribution to the substitution error (when one replaces $T_0\ge \epsilon\mu^{-1}$ by $T_0=\epsilon \mu^{-1}$) of all elements with $\ell_K\ge C_0{\bar\ell}_K$ does not exceed $C\mu^{-1}h^{-3}+ C\mu^2 h^{-2}{\bar\ell}_K^{-1}$ which results in the same estimate.

\medskip
\noindent
(iii) Alternatively, assume that $\ell_j\le \gamma\le \ell_{j+1}$ for some $j\le K-1$. Again one needs to consider elements with $\ell_k \asymp {\bar\ell}_k\le \gamma $ for some $k\le j$. Then in virtue of the same arguments of the proof of proposition \ref{IRO8-prop-4-5} of \cite{IRO8} the contribution to the asymptotics does not exceed 
\begin{equation*}
C\mu^2h^{-2}\ell_k^{-1} \gamma^{-1}|\log h|^J\times \gamma^4\times \ell_k\gamma^{-1} \asymp C\mu^2h^{-2} |\log h|^J\gamma^2
\end{equation*}
where the last factor $\ell_k\gamma^{-1}$ is the upper bound of the relative measure of $\ell_k$ elements to $\gamma$.

Again summation with respect to $\gamma$ results in $C\mu^2h^{-2}|\log h|^J$.

\medskip
\noindent
(iv) Finally, consider elements with $\ell_1\ge \gamma$. Then I redefine $\gamma =\ell_1$ and only $\ell\asymp {\bar\ell}_1$ should be considered since otherwise $T_0\le \epsilon\mu^{-1}$.

Contribution of such elements to the remainder with $T_0=\epsilon\mu^{-1}$   does not exceed $C\mu h^{-3}{\bar\ell}_1^4=O(\mu^3h^{-1}|\log h|^J)$. 
\end{proof}

\subsubsection
\label{sect-4-1-3}
Now I want to attack resonances. 

\begin{proposition}\label{prop-4-4}  Assume that that $|f_1-f_2|\ge \epsilon_0$ and  $|\nabla (f_1f_2^{-1})|\ge \epsilon_0$  on $\supp\psi$. Let consider asymptotics with the main part given by the standard implicit formula with any $T\ge \epsilon \mu^{-1}$. Then under condition $(\ref{4-1})$  
 the remainder does not exceed $(\ref{4-2})_q'$\ \footnote{\label{foot-16} Where $'$ denotes that the power of $|\log h|$ could be larger.} with $q<2$ arbitrarily close to $2$.

\end{proposition}

\begin{proof} [Proof, Part I]\label{pf-4-4-I}
Let us consider first zone $\{|x_1|\ge {\bar\gamma}\}$ and apply to the canonical form the same method as in the proof of proposition \ref{prop-4-3}.
Again I assume that there is just one resonance surface $\Xi=\{x_1=0\}$. 

However, definition of $\ell_k$ involves derivatives of order $k$ and they are unbounded (with respect to $x_1$) as $\gamma^{k+1}\le \mu^{4-2m}$. To avoid this problem I rescale first $B(y,\gamma(y))$ into $B(0,1)$ and then apply this method. However, such rescaling would replace $\rho$ by $\rho \gamma$ and $\ell$ by $\ell\gamma$ and leave $\rho\ell\ge C\mu h|\log h|$ intact. However then the main part of remainder estimate gets factor $\gamma^{-1}$ since the number of indices ``$p$'' will be $\rho(\mu h\gamma)^{-1}$ because the derivative with respect to $x_1$ for ``$p$'' and ``$(p+1)$'' would differ by $\mu h\gamma$.

So, the contribution of zone $\{x: \gamma(x)\asymp \gamma\}$ to the error estimate in question does not exceed $C (\mu h)^qh^{-4}$ with $q<2$ arbitrarily close to $2$; factor $\gamma^{-1}$ discussed above is compensated by the  same factor coming from the measure. Summation with respect to $\gamma$ results in the similar answer (extra factor $|\log h|$ is covered by miniscule decrease of $q$).

This works as long as $|x_1|\ge {\bar\gamma}$ with
\begin{equation}
{\bar\gamma}=\mu^{2-m+\delta'}+\mu^{-2} 
\label{4-5}
\end{equation}
where the last term dominates as $m\ge5$ only and guarantees that $\gamma\ge C\mu^{-1}h^{1-\delta}$.

Therefore,

\begin{claim}\label{4-6}
Contribution of zone $\{{\bar\gamma}\le |x_1|\le \epsilon\}$ to the remainder does not exceed $C\mu^{-1}h^{-3}+C(\mu h)^q h^{-4}$.
\end{claim}

In particular, as $m=5$ proposition is proven since the contribution of zone $\{|x_1|\le {\bar\gamma}$ to the remainder does not exceed $C\mu h^{-3}{\bar\gamma}$.
\end{proof}

\begin{proof} [Proof of Proposition 4.4, Part II]\label{pf-4-4-II}
On the other hand, in zone $\{|x_1|\le {\bar\gamma}\}$ one can use precanonical form and contribution of subzone  $\{\rho \ge {\bar\rho}\Def C\max ({\bar\gamma}, \mu h|\log h| {\bar\gamma}^{-1})\}$ to the error is negligible while contribution of subzone  
$\{\rho \le {\bar\rho}\}$ to the remainder does not exceed $C\mu h^{-3}{\bar\gamma}{\bar\rho}= C\mu h^{-3}{\bar\gamma}^2 + C\mu ^2h^{-2} |\log h|$; picking ${\bar\gamma}$ as $m=4$ results in the proper estimate then.

\medskip
However as $m=3$ one recovers only remainder estimate $C\mu^{-1+2\delta'} h^{-3}+C(\mu h)^qh^{-4}$ which is the required estimate as $\mu\ge h^{-1/3+\delta'}$ and only marginally worse otherwise. To recover proper estimate one needs to reexamine zone $\{|x_1|\le {\bar\gamma},\rho\le {\bar\gamma}\}\setminus \{|x_1|\le {\bar\gamma}_0,\rho\le {\bar\gamma}_0\}$
since the contribution to the remainder of the latter does not exceed $C\mu^{-1}h^{-3}$. Without any loss of the generality one can assume that non-diagonal term does not depend on $x_1$ since one can remove term divisible by $x_1$ by the same method as for $|x_1|\ge {\bar\gamma}$ was removed the whole term.

In subzone $\{\max(\rho,{\bar\gamma}_0)\le |x_1|\le {\bar\gamma}\}$ one can pick up $T_0= Ch|\log h|\rho^{-2}$ which is less than $\epsilon \mu^{-1}$ unless $\rho\le {\bar\rho}\Def C(\mu h|\log h|)^{1/2}$ (which is less than $\mu^{-1}$ and the contribution to the remainder of $\{ |x_1|\le {\bar\gamma},\ \rho \le  {\bar\rho}\}$ does not exceed 
\begin{equation*}
C\mu h^{-4}{\bar\rho}{\bar\gamma}= C\mu^{-1}h^{-3}\times \mu^{\delta}(\mu^3h|\log h|)^{1/2}\le C\mu^{-1}h^{-3}
\end{equation*}
as $\mu\le h^{-1/3+\delta'}$ and $\delta=\delta '$.

In subzones $\{|x_1|\asymp \gamma \le \rho \}$ with 
${\bar\gamma}_0\le \gamma\le {\bar\gamma}_1$ and $\{|x_1|\le \gamma ={\bar\gamma}_0\le \rho\}$ one can take $T_0=Ch|\log h|(\rho\gamma)^{-1}$ which is less than $\epsilon \mu^{-1}$ here for sure. 
\end{proof}

\begin{remark}\label{rem-4-5} It can happen that $T_1$ described in section \ref{sect-3} is less than $\epsilon\mu^{-1}$. Then one can take $T_0\le \epsilon \mu^{-1}$ anyway. 
\end{remark}

\subsubsection{}
\label{sect-4-1-4}
Let us derive sharp remainder estimates in the resonance case:
 
 \begin{proposition}\label{prop-4-6}   Assume that that $|f_1-f_2|\ge \epsilon_0$ and also $|\nabla (f_1f_2^{-1})|\ge \epsilon_0$  on $\supp\psi$. Further, assume that conditions $(\ref{4-1})$ and \ref{3-21-q}  and fulfilled with $q> 2$ and that on (any) resonance surface $\Xi$
 \begin{phantomequation}\label{4-7}\end{phantomequation}
\begin{multline}
 \Hess _\Xi \bigl(\mu h f V^{-1}-1\bigr)  \text{\ \ has at least $r$ eigenvalues  }\\ \text{with absolute values greater than $\epsilon_0$}
\tag*{$(4.7)_r$}\label{4-7-r}
\end{multline}
with $r=q-1$. 

Then the remainder estimate is given by $(\ref{4-2})'_q$   while the main term of asymptotics given by the standard implicit formula with any $T\ge \epsilon \mu^{-1}$.
\end{proposition}

\begin{proof} Due to proposition \ref{prop-4-3} one should cover only case $\mu \ge h^{-1/3+\delta'}$.  Assumptions of proposition imply that condition $(\ref{3-66})_{r+1,r}$ holds. While one can apply proposition \ref{prop-3-10} directly only as $r> 2$, the proof of it yields that  under condition $(\ref{4-7})$ with $r=1,2$ the remainder estimate $O\bigl(\mu^{-1}h^{-3}+ \mu^{2-r/2}h^{2+r/2}|\log h|^J\bigr)$  holds, which is not worse than $(\ref{4-2})_{r+1}$.
Let us apply the same partition to zones as in propositions \ref{prop-3-9}, \ref{3-10}.

Then zone  
\begin{equation*}
\{\gamma \ge \max(\rho, {\bar\gamma}_2\Def C\max(\mu^{-2/3},(\mu h|\log h|)^{1/2}))\}
\end{equation*}
is covered by the arguments used in  the proof of proposition \ref{prop-4-1}; its contribution to the error does not exceed \ref{4-2-q}. Meanwhile in the subzone  $\{\rho \ge \gamma \ge  {\bar\gamma}_2)\}$ one can take $T_0= Ch|\log h|{\bar\gamma}_2^{-2}$ which is less than $\epsilon\mu^{-1}$.

Further, in subzone $\{{\bar\gamma}\le |x_1|\asymp \gamma \le {\bar\gamma}_2\}$ one can take $T_0= Ch|\log h|\zeta^{-2}$ which is less than $\epsilon\mu^{-1}$ unless $\zeta \le C(\mu h|\log h|)^{1/2}$. Similarly, in the subzone $\{|x_1|\le {\bar\gamma}\}$ one can take $T_0= Ch|\log h|\zeta^{-2}$ as well. Note that the contribution of the zone 
\begin{equation*}
\{|x_1|\le {\bar\gamma}_2, \zeta \le {\bar\zeta}\Def C(\mu h|\log h|)^{1/2}\}
\end{equation*}
to the remainder does not exceed $C\mu h^{-3}{\bar\gamma}_2 {\bar\zeta}^r$ which does not exceed the second term in $(\ref{4-2})_{r+1}$ as ${\bar\zeta}\asymp {\bar\gamma}_2$ i.e. $\mu^{-2/3}\le C(\mu h|\log h|)^{1/2}$ i.e. as $\mu \ge C(h|\log h|)^{-3/7}$. 

As $\mu \le C(h|\log h|)^{-3/7}$ let us consider zone 
$\{|x_1|\le {\bar\gamma}_2\}$ and use the precanonical form here. Then 
contribution of the subzone $\{|x_1|\le {\bar\gamma}_2, \ \rho \le {\bar\rho}_2, \ \zeta\le {\bar\zeta} \}$ to the remainder does not exceed 
$C\mu h^{-3}{\bar\gamma}_2 {\bar\rho}_2{\bar\zeta}^r$ which does not exceed the second term in $(\ref{4-2})_{r+1}$ as ${\bar\zeta}\ge {\bar\gamma}_2{\bar\rho}_2$. Therefore  one can pick up ${\bar\rho}_2={\bar\zeta}{\bar\gamma}_2^{-1}$; one can see easily that ${\bar\rho}_2\gg \mu^{-2/3}$ and ${\bar\rho}_2{\bar\gamma}_2\gg C\mu h|\log h|$ and then $T_0\le Ch|\log h|{\bar\rho}_2^{-1}{\bar\gamma}_2^{-1}\ll \epsilon \mu^{-1}$ and the contribution of the zones
$\{|x_1|\le {\bar\gamma}_2, \ \rho \ge {\bar\rho}_2 \}$ and
$\{|x_1|\le {\bar\gamma}_2,  \ \zeta\ge {\bar\zeta} \}$ to the error is negligible.
\end{proof}

\subsubsection{}
\label{sect-4-1-5} Finally, let us consider the vicinity of $\Sigma$.

\begin{proposition}\label{prop-4-7} Let condition $(\ref{1-1})-(\ref{1-2}),(\ref{4-1})$ be fulfilled and let $\psi$ be  supported in the small vicinity of $\Sigma$.

Then the remainder estimate is given by $(\ref{4-2})'_q$ with  any $q<2$ arbitrarily close to $2$ while the main term of asymptotics given by the standard implicit formula with any $T\ge \epsilon \mu^{-1}$.
\end{proposition}

\begin{proof} Again this remainder estimate is no smaller than $O(\mu^{-1}h^{-3}+\mu^2h^{-2})$ delivered by proposition \ref{prop-3-13}. 

Applying the same approach as in the Part I of the proof of  proposition \ref{prop-4-4} (which is applicable as $\gamma \ge {\bar\gamma}_2\Def \max \bigl((\mu h|\log h|)^{1/2}, \mu^{-1}\bigr)$ now) I conclude after summation with respect to $\gamma$ that the contribution of the zone $\{x:\gamma (x)\asymp \gamma\}$ (where $\gamma(x)=\dist (x,\Sigma)\}$) to the error does not exceed $C(\mu h)^q h^{-4}\gamma$ with  $q<2$ arbitrarily close to $1$. Here an extra factor $\gamma$ appears because the measure of $\{x:\gamma(x)\asymp\gamma\}$ is $\asymp\gamma^2$ rather than $\asymp\gamma$ as it was before. Then after summation I conclude that

\begin{claim}\label{4-8}
The contribution of the zone $\{x: \max ((\mu h, \mu^{-1})\le \gamma(x)\le \gamma\}$ to the error is $O\bigl((\mu h)^q h^{-4}\gamma\bigr)$. In particular, it is $O\bigl((\mu h)^q h^{-4}\bigr)$ as $\gamma =\epsilon$.
\end{claim}

Meanwhile, the contribution of the zone $\{x: \gamma(x)\le {\bar\gamma}_2\}$ to the error does not exceed $C\mu h^{-3} {\bar\gamma}_2^2\asymp (\mu^{-1}h^{-3}+\mu^2h^{-2})$.
\end{proof}

\subsubsection{}
\label{sect-4-1-6}
Let us improve the above estimate under generic conditions.

\begin{proposition}\label{prop-4-8} Let conditions $(\ref{1-1})-(\ref{1-2}),(\ref{4-1})$ and 
\ref{3-88-r} with $r=1,2$ be fulfilled. Let $\psi$ be supported in the small enough vicinity of $\Sigma$. Then

\medskip
\noindent 
{\rm (i)} The total  remainder is given by $(\ref{4-2})'_{r+1}$ while the main part is given by the standard implicit formula with any $T\ge \epsilon \mu^{-1}$. 

\medskip
\noindent 
{\rm (ii)} Under condition $(\ref{3-85})$ the total  remainder is given by $(\ref{4-2})'_{r+2}$ as $r=0,1$ while the main part is given by the standard implicit formula with any $T\ge \epsilon \mu^{-1}$;

\medskip
\noindent 
{\rm (iii)} Under conditions $(\ref{3-85})$ and $(\ref{3-84})$ the total  remainder is given by $(\ref{4-2})'_{r+2}$  as $r=2$ while the main part is given by the standard implicit formula with any $T\ge \epsilon \mu^{-1}$.
\end{proposition}

\begin{proof} Again let us note that the remainder estimate given by proposition \ref{prop-3-17} as $r=0$ and  \ref{prop-3-18} as $r=1,2$ is no worse than one announced here. 

In zone $\{\gamma \ge \max (\rho,(\mu h|\log h)^{1/2},\mu^{-1})\}$ the same arguments as in the proof of proposition \ref{prop-4-3} are applied, and conditions \ref{3-88-r}  and (\ref{3-85}) add factors $ (\mu h|\log h|)^{r/2}$ and  $ (\mu h|\log h|)^{1/2}\gamma^{-1}$ respectively to the measure of zone $\{\rho\le (\mu h|\log h|)^{1/2}$.  

Then after summation with respect to $\gamma$ one arrives to the estimate $(\ref{4-2})'_{r+k+1}$ of the contribution of this zone to the error estimate and $k=1$ under condition (\ref{3-85}) and $k=0$ otherwise. 

Meanwhile $T_0\le \epsilon\mu^{-1}$ in zone 
$\{\rho\ge \gamma \ge C\max ((\mu h|\log h|)^{1/2} ,\mu^{-1})\}$ .

So, one needs to consider zone $\{\gamma\le C{\bar\gamma}_2\}$. In this zone condition (\ref{3-75}) has no value. Contribution of this zone to the remainder estimate obviously does not exceed $(\ref{4-2})'_2$. Therefore case $r+k\le 1$ is covered.
Let us introduce $\zeta$ as before. As \emph{either} $\zeta\le \gamma$ \emph{or} condition (\ref{3-84}) is fulfilled, one can take $T_0= Ch\zeta^{-2}$ and then only subzone $\{\zeta \le C(\mu h|\log h|)^{1/2}\}$ should be considered; then condition \ref{3-88-r} adds an extra factor $(\mu h|\log h|)^{r/2}$ to the measure and to the estimate which becomes $(\ref{4-2})'_{r+2}$.

That leaves us with the analysis of zone $\{\zeta\ge \gamma\}\cap \{\gamma\le {\bar\gamma}_2\}$ and only without condition (\ref{3-84}), in which case only estimate $(\ref{4-2})'_4$ should be proven under condition $(\ref{3-88})_1$. However then one can take $T_0=C\mu h|\log h|\gamma^{-1}\zeta^{-1}$ and then only subzone $\zeta\le C\mu h|\log h|\gamma^{-1}$ remains to be considered. Its contribution to the remainder does not exceed $C\mu h^{-3}\times \mu h|\log h|\gamma^{-1}\times \gamma^2$ which sums with respect to $\gamma$ ranging from 
${\bar\gamma}_1\Def \max(C\mu h|\log h|,\mu^{-1})$ to ${\bar\gamma}_2$ to $C\mu^2 h^{-2}|\log h|{\bar\gamma}_2$ which is properly estimated. Meanwhile contribution of $\{x:\gamma(x)\le {\bar\gamma}_1\}$ to the remainder does not exceed $C\mu h^{-3}{\bar\gamma}_1^2$ and is properly estimated as well.
\end{proof}

\subsubsection{Conclusion}
\label{sect-4-1-7}

\begin{corollary}\label{cor-4-9}  Let conditions of  one of propositions \ref{prop-4-3} -- \ref{prop-4-8} be fulfilled. Further, let condition $(\ref{0-8})$  be fulfilled. Then 

\medskip
\noindent
{\rm (i)} The remainder estimate is $(\ref{4-2})'_q$ where 
\begin{enumerate}[label=\rm (\alph*)]
\item $q<2$ is arbitrarily close to $2$ in the general case (propositions \ref{prop-4-3}, \ref{prop-4-4}, \ref{prop-4-7}) and 
\item $q$ is described in the corresponding case (propositions \ref{prop-4-1}, \ref{prop-4-6}, \ref{prop-4-8}), 
\end{enumerate}
while the main part of asymptotics is given by $(\ref{4-4})$.

\medskip
\noindent
{\rm (ii)} In particular  the remainder estimate is $O(\mu^{-1}h^{-3})$ as $\mu \le {\bar\mu}_q\Def Ch^{-q/(q+4)}|\log h|^{-K}$. In particular 
\begin{enumerate}[label=\rm (\alph*)]
\item in the general case ${\bar\mu}_1= Ch^{\delta-1/3}$ 
\item in the generic case ${\bar\mu}_4= h^{-1/2}|\log h|^{-K}$.
\end{enumerate}
\end{corollary}

\begin{proof} Proof coincides with the proof of corollary \ref{cor-4-2}. \end{proof}

Then I get immediately 

\begin{corollary}\label{cor-4-10}  {\rm (i)} Theorem \ref{thm-0-2} is proven with $\cE^\MW_\corr=0$ as $\mu \le h^{\delta-1/3}$;

\medskip
\noindent
{\rm (ii)} Theorem \ref{thm-0-3} is proven  as $\mu \le Ch^{-1/2}|\log h|^{-K}$.
\end{corollary} 

\subsection{Strong magnetic field}
\label{sect-4-2}

In this subsection I assume that
\begin{equation}
h^{-1/3}\le \mu \le ch^{-1}
\label{4-9}
\end{equation}
sometimes making separate considerations for the case of the \emph{superstrong\/} magnetic field  
\begin{equation}
h^{-1}|\log h|^{-K}\le  \mu \le ch^{-1}
\label{4-10}
\end{equation}
as needed.   Also in the general case I consider a special range $h^{-1/3+\delta}\le \mu \le h^{-1/3-\delta}$.

\subsubsection{}
\label{sect-4-2-1}
I start from the regular points. 

\begin{proposition}\label{prop-4-11}  Assume that  there are no resonances of order not exceeding $M$ and also $|\nabla (f_1f_2^{-1})|\ge \epsilon_0$  on $\supp\psi$. Further, assume that condition \ref{3-21-q} with $q\ge 1$ is
fulfilled and
\begin{equation}
\mu^{-1}h|\log h|\le \varepsilon \le \mu h
\label{4-11}%+2
\end{equation}
with $\varepsilon = \mu^{-2}$ here.

Then  the remainder does not exceed
\begin{phantomequation}\label{4-12}\end{phantomequation}
\begin{equation}
C\mu^{-1}h^{-3}+  Ch^{-4}(\mu h|\log h|) ^{q/2} \varepsilon+ 
\left\{\begin{aligned} 
&0 \qquad &&q\ge 3,\\ 
&\mu^2h^{-2}\varepsilon^{(q-1)/2} &&2\le q<3,\\
&\bigl(\mu^2h^{-2}\varepsilon^{(q-1)/2} + \mu h^{-3}\varepsilon^{q/2}\bigr)\qquad&&1\le q<2
\end{aligned}
\right.
\tag*{$(4.12)_q$}
\label{4-12-q}
\end{equation}
 while the main term of asymptotics given by $(\ref{4-4})$.
\end{proposition} 
\begin{proof} 
The proof follows the sequence of the proofs of propositions \ref{IRO8-prop-4-15}, \ref{IRO8-prop-4-17} with $\varepsilon=\mu^{-2}$, 
\ref{IRO8-prop-4-20}(i) from \cite{IRO8}: 

\medskip
\noindent
(i) First, using the method of successive approximations, I rewrite the implicit formula as expression (\ref{4-4}) plus a (temporary) \emph{correction term\/}
\begin{equation}
\sum_{\iota} 
\int \Bigl({\widehat\cE}^\MW_{\bar Q_\iota} (x,0)- \cE^\MW_{\bar Q_\iota} (x,0)\Bigr)\,dx
\label{4-13}%+3
\end{equation}
with an integrand defined by (\ref{IRO8-4-33})-(\ref{IRO8-4-34}), \cite{IRO8}:
\begin{equation}
{\widehat\cE}^\MW_{\bar Q} (x,0) \Def  \const \sum_{n,p} \Bigl(\theta (\cA_{pn}){\bar Q}\Bigr)\Bigr|_{(x'',\xi'')=\Psi^{-1}(x)} \times 
f_1(x)f_2(x) \mu^2h^{-2}
\label{4-14}
\end{equation}
and
\begin{multline}
\cE^\MW_{\bar Q}(x,0)\Def \\\const \sum_{n,p} \theta\Bigl(V(x)- (2n+1)\mu hf_2(x)- (2p+1)\mu h f_1(x)\Bigr) f_1(x)f_2(x){\bar Q} \mu^2 h^{-2}
\label{4-15}
\end{multline}
where $Q_\iota$ cover zone with $T_0\ge \epsilon\mu^{-1}$ (i.e. $\rho\le C(\mu h|\log h|)^{1/2}$); here $\const$ meant the same constant as in the definition of $\cE^\MW$.

\medskip
\noindent
(ii) Then I need to estimate  expression (\ref{4-13}). Now, however, one needs  only analysis which was used in \cite{IRO8} in the strictly outer zone. The crucial moment is the estimate of the correction term as in (\ref{IRO8-prop-4-17}) of \cite{IRO8} which was
\begin{multline}
C\mu^2h^{-2} \times h\rho^{-2}\times \varepsilon h^{-1}\times 
\bigl(\rho^2(\mu h)^{-1}+1\bigr)\times \bigl(\rho(\mu h)^{-1}+1\bigr)\rho^{q-1}\asymp \\
C\mu^2h^{-2} \varepsilon \bigl(\rho^2(\mu h)^{-1}+1\bigr)\times \bigl(\rho(\mu h)^{-1}+1\bigr)\rho^{q-3}
\label{4-16}
\end{multline}
and one needs to sum this expression with respect to $\rho$ ranging from $\varepsilon^{1/2}$ ($=\mu^{-1/2}$  as $\varepsilon=\mu^{-2}$)  to ${\bar\rho}\Def C(\mu h|\log h|)^{1/2}$. Then 

\begin{enumerate}[label=(\alph*)] 

\item
As $q>3$  expression (\ref{4-16}) sums to its value as $\rho={\bar\rho}$, which is exactly \newline$C\varepsilon (\mu h|\log h|)^{q/2}h^{-4}$. 

\item As $1<  q\le 3$ an extra term due to summation of $C\mu ^2 h^{-2}\varepsilon \rho ^{q-3}$ appears in the estimate.  As $q<3$ this results in the value of this term as $\rho=\varepsilon^{1/2}$ which is $C\mu ^2 h^{-2}\varepsilon ^{(q-1)/2}$ while as $q=3$ it results in 
$C\mu^2h^{-2}\varepsilon \log (\mu h |\log h| \varepsilon^{-1})$.

\item Further, as $1\le q\le 2$ one should take in account also 
$C\mu h^{-3}\varepsilon  \rho ^{q-2}$ which sums to 
$C\mu h^{-3}\varepsilon  ^{q/2}$ as $1<q<2$ and $C\mu h^{-3}\varepsilon  |\log h|$ as $q=2$.
\end{enumerate}

\noindent
(iii) Finally  in zone $\{\rho \le \varepsilon^{1/2}\}$ one simply considers its contribution to the remainder with $T=\epsilon \mu^{-1}$ rather than to the correction:
\begin{multline}
C\mu^2h^{-2} \times 
\bigl(\rho^2(\mu h)^{-1}+1\bigr)\times \bigl(\rho(\mu h)^{-1}+1\bigr)\rho^{q-1}\bigr|_{\rho=\varepsilon^{1/2}}\asymp \\
C\mu h^{-3} \varepsilon^{q/2}+C\mu^2h^{-2} \varepsilon^{(q-1)/2}.
\label{4-17}
\end{multline}
\end{proof}

To cover properly the case of the \emph{superstrong\/} magnetic field  (\ref{4-10}) let us  rewrite the implicit formula as
\begin{equation} 
\int {\widehat\cE}^\MW_I (x,0)\,dx
\label{4-18}
\end{equation}
and then as $q>3$ rewrite it as the three term decomposition with respect to $\varepsilon$; then the third (remainder) term is $o(\mu^{-1}h^{-3})$ while the second term is given by a two-dimensional Riemannian sum proportional $\mu^{-2}$ with the steps $2f_1\mu h$ and $2f_2\mu h$.  Replacing this Riemannian sum by an integral one can see easily that again with an error $o(\mu^{-1}h^{-3})$ the second term is $\varkappa \mu^{-2}h^{-4}$ which disagrees with expression as $\mu \le h^{-1}|\log h|^{-K}$ unless $\varkappa=0$. Combining with the previous proposition I arrived to

\begin{proposition}\label{prop-4-12}
In frames of proposition \ref{prop-4-11} with $q>3$ asymptotics with the main part $(\ref{4-4})$ and $O(\mu^{-1}h^{-3})$ remainder holds.
\end{proposition}

Combining  propositions \ref{prop-4-11}, \ref{prop-4-12}, corollary \ref{cor-4-2} and results of section \ref{sect-2}  I conclude that 

\begin{claim}\label{4-19}
in frames of proposition \ref{prop-4-11} in the generic case $q>3$ estimate $O(\mu^{-1}h^{-3})$ is proven for the complete range of $\mu$ ($1\le \mu \le ch^{-1}$).  Meanwhile in the general case $q=1$ estimate $O(\mu^2h^{-2})$ is proven as $h^{-1/2}|\log h|^K\le \mu \le ch^{-1}$. 
\end{claim}

\subsubsection{}
\label{sect-4-2-2}
Now I want to improve result in the general case.

\begin{proposition}\label{prop-4-13}  Assume that $f_1\ne f_2$ and there are no resonances of order not exceeding $M$  and also critical points of $f_1f_2^{-1}$  are non-degenerate on $\supp\psi$. Then as $\mu \ge h^{-1/3}$ the remainder estimate is given by \ref{4-12-q} with  any $q<2$ arbitrarily close to $2$ while the main term of asymptotics given by  $(\ref{4-4})$.
\end{proposition}

\begin{proof} 
Proof follows ideas of the proof of proposition \ref{IRO8-prop-4-22} of \cite{IRO8}. Let us introduce functions $\ell_k$, $k=1,\dots,K$  in the same proof but with $\gamma=1$ and again let ${\bar\ell}_k= (\mu h|\log h|^J)^{1/(k+1)}$. 

\medskip
\noindent
(i) Assume first that $|\nabla (f_1f_2^{-1})|\ge \epsilon_0$. Then  the arguments of the mentioned proof survive with this simplification. 

\medskip
\noindent
(ii) Assume now that $f_1f_2^{-1}$ has one non-degenerate stationary point ${\bar x}$.   Assume first that $\ell_K\le \gamma$ where $\gamma ={\frac 1 2}|x-{\bar x}|$. Then the number of indices ``$p$'' does not exceed $C\ell_k^k/(\mu h \gamma)$ and the contribution of all such elements with $\ell_k\asymp {\bar\ell}_k$ for some $k=1,\dots,(K-1)$ to the correction term does not exceed $C(\mu h)^{1-\delta}\cdot \mu^{-2} h^{-4}\gamma^{-1}\times \gamma^4|\log h|^J$ where $\gamma^4$ is the total measure of such elements, $\delta=1/K$ and $J$ is large enough. Summation with respect to $\gamma$ results in the announced estimate.

On the other hand, the total contribution to the correction of all elements with $\ell_k\ge C_0{\bar\ell}_k$ does not exceed $C\mu^{-1}h^{-3}+ C\mu  h^{-3}\cdot \mu^{-2}{\bar\ell}_K^{-1}$.

\medskip
\noindent
(iii) Alternatively, assume that $\ell_j\le \gamma\le \ell_{j+1}$. Again one needs to consider elements with $\ell_k \asymp {\bar\ell}_k\le \gamma $ for some $k\le j$. Then the contribution to the correction term does not exceed $C\mu h^{-3} \times \mu^{-2} \gamma^{-1}\times \gamma^4\times \gamma^{-1} |\log h|$ where the second factor $\gamma^{-1}$ appears since the relative measure of $\ell_k$ elements to $\gamma$ elements does not exceed $\ell_k\gamma^{-1}$. Again summation with respect to $\gamma$ results in $C\mu ^{-1} h^{-3}|\log h|^J$.

\medskip
\noindent
(iv) Finally, consider elements with $\ell_1\ge \gamma$. Then I redefine $\gamma =\ell_1$ and only $\ell_1\asymp {\bar\ell}_1$ should be considered. Contribution of such elements to the correction does not exceed $C h^{-4}\times\mu^{-2}{\bar\ell}_1^4=O(\mu^{-1}h^{-3})$. 
\end{proof}

Therefore

\begin{claim}\label{4-20}%+9
In the general case  asymptotics with the main part (\ref{4-4}) and the remainder estimate $O(\mu^{-1}h^{-3}+\mu^2h^{-2})$ holds unless 
\begin{equation}
h^{-1/3+\delta}\le \mu \le h^{-1/3-\delta}
\label{4-21}
\end{equation}
in which case the remainder estimate is $O(\mu^{-1}h^{-3-\delta})$ containing an extra factor $h^{-\delta}$; $\delta>0$ is an arbitrarily small exponent. 
\end{claim}

\subsubsection{}
\label{sect-4-2-3}
Let us  recover  remainder estimate $O(\mu^{-1}h^{-3}+\mu^2h^{-2})$ in the latter case (\ref{4-21}) (introducing some correction term). Using the same arguments as above I can purge from $\cA_{pn}$ in (\ref{4-14}) all terms which are even marginally less than $\mu^{-2}$; it includes higher order terms and also the difference between $\mu^{-2}\cB_{pn}$ and $\mu^{-2}\cB_{{\bar p},{\bar n}}$ where ${\bar p}=\alpha /(2f_1\mu h)$ and
${\bar n}=(1-\alpha)/(2f_2\mu h)$, $\alpha=\alpha(x)$ is the minimizer of $|\nabla \phi_\alpha|^2$.

Then $\cB_{pn}$ becomes $\cB_{{\bar p},{\bar n}}=\omega (x)$\,\footnote{\label{foot-17} Which is the smooth function outside of the critical points of $f_1f_2^{-1}$.}. However, let us include this modified term $\mu^{-2}\cB_{pn}$ for all $p$, $n$ and not only for those for which $T_0\ge \epsilon \mu^{-1}$. Then one needs to correct this alternation  by the term $\varkappa \mu^{-2}h^{-4}$ with $\varkappa$ selected so it would result in the correction term 0, if one replaces the Riemann sum by the corresponding integral because it would provide the result for $\mu \ge h^{-1/3-\delta}$ and it should agree with the results of the previous subsubsection. I leave the easy details to the reader.  Then I arrive to

\begin{proposition}\label{prop-4-14}  Assume that $f_1\ne f_2$ and there are no resonances of order not exceeding $M$  and also the critical points of $f_1f_2^{-1}$  are non-degenerate on $\supp\psi$. Then for $\mu$ satisfying $(\ref{4-21})$  the remainder estimate is $O(\mu^{-1}h^{-3}+\mu^2h^{-2})$  while the main part of the asymptotics is given by  
\begin{equation}
\int \Bigl(\cE^\MW(x,0)+ \cE_\corr^\MW(x,0)\Bigr) \psi (x)\, dx.
\label{4-22}
\end{equation}
with
\begin{align}
\cE^\MW_\corr (x,\tau)=\qquad\qquad\ & \label{4-23}\\
(2\pi)^{-2}\mu^2h^{-2}\sum _{(p,n)\in \bZ^{+\,2} } 
\Bigr(&\theta \bigl(2\tau+V- (2p+1)\mu h f_1 - (2n+1)\mu hf_2-\omega (x)\mu^{-2}\bigr) -\notag\\
&\theta \bigl(2\tau+V- (2p+1)\mu h f_1 - (2n+1)\mu hf_2\Bigr) f_1 f_2 \sqrt g+\notag \\
&\qquad\qquad\qquad\qquad\qquad\qquad(4\pi)^{-2}\mu^{-2}h^{-4}(2\tau +V)\omega \sqrt g\notag.
\end{align}
\end{proposition}

\begin{remark}\label{rem-4-15} (i) Here 
\begin{equation*}
\int \cE^\MW_\corr (x,\tau)\psi(x)\, dx=O(\mu^{-1}h^{-3-\kappa})
\end{equation*}
with an arbitrarily small exponent $\kappa>0$\,\footnote{\label{foot-18} It follows from proposition \ref{prop-4-3}};

\medskip
\noindent
(ii)  Furthermore under nondegeneracy condition 
\begin{equation}
\sum_{|\beta|\le K}|\nabla^\beta \phi_\alpha|\ge \epsilon_0\qquad \forall x\ \forall \alpha\in[0,1]
\label{4-24}
\end{equation}
with arbitrarily large $K$ one can skip a correction term  without deteriorating remainder estimate unless $h^{-1/3}|\log h|^{-J}\le \mu \le h^{-1/3}|\log h|^J$ and with the remainder estimate $O(h^{-8/3}|\log h|^J)$ in this border case
\footnote{\label{foot-19} Notice   that  under  condition (\ref{4-24}) ${\bar\ell}_K\asymp 1$ in the proof  of propositions \ref{prop-4-3}, \ref{prop-4-13}.}. I suspect that one can get rid off logarithmic factors and to prove an  estimate $O(\mu^{-1}h^{-3}+\mu^2h^{-2})$ even in the border case.

\medskip
\noindent
(iii) However I could not find any example demonstrating that this correction term is not superficial and without (\ref{4-24}) one cannot skip it without penalty.  Clarification of this would be interesting.

Correction term in \cite{IRO8} was not superficial for sure. 
\end{remark}

\subsubsection
\label{sect-4-2-4}
Now I want to attack resonances. I start from the generic case. First, after rescaling again in the same manner as before one can see easily that the contribution of zone $\{|x_1\asymp\gamma\}$ to an error does not exceed $C\mu h^{1-d}\varepsilon^{1-\kappa}$ where factor $\gamma^{-1}$ appearing from the calculation of the number of indices ``$p$'' is compensated by factor $\gamma$ appearing from the measure. Then this contribution does not exceed $C\mu h^{-3}(\mu^{4-2m}\gamma^{-1})^{(1-\kappa)}$ if 
only terms originated from non-diagonal terms are removed and summation with respect to $\gamma\ge {\bar\gamma}$ results in the value of this expression as $\gamma={\bar\gamma}$ and as $m\ge 4$ it is $O(\mu^{-1}h^{-3}+\mu^2h^{-2})$. I remind that ${\bar\gamma}=\mu^{\delta-2}$ as $m=4$ and ${\bar\gamma}=C\mu^{-2}$ as $m\ge 5$. 

On the other hand in the zone $\{|x_1|\le {\bar\gamma}\}$ one can apply proposition \ref{prop-3-10} and get $C(\mu h)^{3/2}h^{-4}{\bar\gamma}|\log h|^J$ which is $o(\mu^{-1}h^{-3}+\mu^2h^{-2})$ as well (as $m\ge 4$). 

Now after singular terms from operator are removed it can be treated as if there was no resonance resulting in two following statements:

\begin{proposition}\label{prop-4-16}  Assume that $f_1\ne f_2$  and  $|\nabla(f_1f_2^{-1})|\ge \epsilon_0$   on $\supp\psi$. Furthermore, assume that there are no third-order resonances on $\supp\psi$.

Then 

\medskip
\noindent
{\rm (i)} As $\mu \ge h^{-1/3-\delta}$ the remainder estimate is $O(\mu^2h^{-2})$ while the main term of asymptotics given by  $(\ref{4-4})$;

\medskip
\noindent
{\rm (ii)} As $h^{-1/3+\delta}\le \mu \le h^{-1/3-\delta}$ statement (i) of remark \ref{rem-4-15} holds; further under condition $(\ref{4-24})$ one can skip correction term with $O(\mu^2h^{-2}|\log h|^J)$ penalty\footnote{\label{foot-20} Which is probably superficial.}.
\end{proposition}

\subsubsection
\label{sect-4-2-5}
Consider resonances of order 3 now. As $\mu\ge h^{-1/3-\delta}$ in the virtue of the above arguments, the contribution of zone 
\begin{equation}
\bigl\{|x_1|\ge {\bar\gamma}_3\Def \mu^{-3}h^{-1-\delta}\bigr\}
\label{4-25}
\end{equation}
to the correction does not exceed $C\mu^2h^{-2}$.

Moreover, contribution of zone 
\begin{equation}
\bigl\{|x_1|\le {\bar\gamma}_2 \Def (\mu h)^{1/2}|\log h|^{-1/2}\bigr\}
\label{4-26}
\end{equation}
to the correction also
does not exceed $C\mu^2h^{-2}$. However, zone
\begin{equation}
\bigl\{{\bar\gamma}_2 \le |x_1|\le {\bar\gamma}_3\bigr\}
\label{4-27}
\end{equation}
needs to be reexamined; here ${\bar\gamma}_2\ge {\bar\gamma}$ as $\mu\ge h^{-1/3-\delta}$ and as $\mu \ge h^{-3/7-\delta'}$ zone (\ref{4-27}) disappears. In virtue of arguments of subsection \ref{4-2} the contribution of this zone does not exceed $C\mu^2h^{-2-\kappa}$ with arbitrarily small $\kappa>0$ and now I want to improve it marginally.

I claim that

\begin{proposition}\label{prop-4-17}
For any $\delta>0$ there exists $K=K(\delta)$ such that if  $\mu^{-1+\delta}\le \gamma \le \mu^{-\delta}$ and $B(y,\gamma(y))$ with $\gamma(y)\asymp \gamma$ is rescaled to $B(0,1)$ then either
\begin{align}
&\sum_{|\alpha|\le K}  |\partial ^\alpha _z\cA_{pn}|\asymp \varsigma 
\qquad \forall z\in B(0,1)
\label{4-28}\\
\intertext{with some $\varsigma$ or}
&\sum_{|\alpha|\le K}  |\partial ^\alpha _z\cA_{pn}|\le \varsigma 
\qquad \forall z\in B(0,1)
\label{4-29}
\end{align}
with  $\varsigma=\mu h^{1-\kappa}\gamma $ where (continuous) parameters $p$ and $n$ are selected so 
\begin{equation}
\cA_{pn}(y)=\partial_{x_1}\cA_{pn}(y)=0.
\label{4-30}
\end{equation}
\end{proposition}

\begin{proof}
Note that before rescaling $\cA_{pn}=\cA_{pn}^0 + x_1^{-1}\mu^{-2}\cB_{pn} +\dots$ with smooth symbols. Then decomposing all smooth symbols into Taylor series at $y$ one can prove proposition easily since high powers of there contain high powers of $\gamma$ and thus are small.
\end{proof}

\begin{proposition}\label{prop-4-18} In frames of proposition \ref{prop-4-17} contribution of zone $\{\mu^{-1+\delta}\le |x_1|\le \mu^{-\delta}\}$ to the error estimate estimate is $O(\mu^{-1}h^{-3}+\mu^2h^{-2}|\log h|^J)$.
\end{proposition}

\begin{proof} Proof repeats those as I had before. One needs to construct $\ell_k$ corresponding to $\varsigma \cA_{pn}$ and then $(\mu h|\log h|)$ as minimal value for $\ell_k^{k+1}$ should be replaced by  $\varsigma^{-1}(\mu h|\log h|)$ and there will be no final division by ${\bar\ell}_K$ which would be $\asymp 1$. That will give $C\mu^2h^{-2}\times |\log h| \times |\log h|\gamma^{-1}\times \gamma$ where the factors $\gamma^{-1}$ and $\gamma$ appear from division by $(\mu h\gamma)$ and the measure. 

Then summation  with respect to $\gamma$ results in an extra $|\log h|$ factor.
\end{proof}

\begin{remark}\label{rem-4-19} (i)   Under condition 
\begin{equation}
\sum_{|\beta|\le K}|\nabla_\Xi^\beta Vf_1^{-1}|\ge \epsilon_0\qquad \forall x\ 
\label{4-31}
\end{equation}
with arbitrarily large $K$ one can skip a correction term  without deteriorating remainder estimate unless $h^{-1/3}|\log h|^{-J}\le \mu \le h^{-1/3}|\log h|^J$ and with the remainder estimate $O(h^{-8/3}|\log h|^J)$ in this border case.
I suspect that one can get rid off logarithmic factors and to prove an  estimate $O(\mu^{-1}h^{-3}+\mu^2h^{-2})$ even in the border case.

\medskip
\noindent
(ii) One can easily construct a correction term in the case of third-order resonances but an expression seems  rather too complicated. So I leave it to the curious reader.
\end{remark}

\subsubsection{}\label{sect-4-2-6} Let us consider the generic case now:

\begin{proposition}\label{prop-4-20}  Assume that that $|f_1-f_2|\ge \epsilon_0$ and also $|\nabla (f_1f_2^{-1})|\ge \epsilon_0$  on $\supp\psi$. Further let us assume that conditions $(\ref{3-21})_4$ and  $(\ref{4-7})_3$ are fulfilled. 

Then  the remainder is $O(\mu^{-1}h^{-3})$ while the main term of asymptotics given by $(\ref{4-4})$.
\end{proposition}

\begin{proof} Again it is sufficient to consider the case of the single resonance surface $\Xi=\{x_1=0\}$.

\medskip
\noindent
(i) Combining in zone $\{|x_1|\asymp \gamma \ge {\bar\gamma}\Def \mu^{\delta-1}\}$ arguments of the proofs of proposition \ref{prop-4-6} and \ref{prop-4-11} one can estimate contribution of it to the correction  terms by 
$C(\mu h|\log h|)^{q/2}h^{-4}\gamma \times \mu^{-2}\gamma^{-1}$ which after summation with respect to $\gamma$ results in $C(\mu h |\log h|)^{q/2}\mu^{-2} h^{-4}|\log h|$, which in turn is $O(\mu^{-1}h^{-3})$ as long as $\mu\le h^{-1}|\log h|^{-J}$.

\medskip
\noindent
(ii) In zone $\{|x_1|\le {\bar\gamma}\}$ one can use arguments of the proofs of proposition \ref{prop-4-7} and estimate the contribution of this zone to the correction term arising as $T\ge T_0$  is replaced by $T=\epsilon\mu^{-1}$ by $C(\mu h|\log h|)^{5/2}h^{-4}{\bar\gamma}$ which is $O(\mu^{-1}h^{-3})$ as $\mu \le h^{-3/5+\delta'}$. 

\medskip
\noindent
(iii) Let $\mu \ge h^{-5/3+\delta'}$. In zone $\{|x_1|\le {\bar\gamma}\}$ one can consider precanonical form and then the contribution of this zone to the correction term does not exceed 
$C(\mu h|\log h|)^{3/2}\mu^{-1}h^{-4}{\bar\gamma} $ which is $O(\mu^{-1}h^{-3})$ as $\mu \le h^{-1+\delta'}$. 

The direct calculation shows that the first approximation term actually vanishes as it comes from the ``main perturbation term'' and is odd with respect to $x_1$, and its estimate contains an extra factor ${\bar\gamma}$ otherwise.

So the correction terms associated with zone $\{|x_1|\le {\bar\gamma}\}$ do not exceed 
\begin{equation*}
C(\mu h)^{3/2}\mu^{-1}h^{-4}{\bar\gamma}^2+ 
C(\mu h)^{1/2}\mu^{-2}h^{-4}{\bar\gamma}=o(\mu^{-1}h^{-3})
\end{equation*}

\medskip
\noindent
(iv) Finally, arguments of (i) should be slightly improved as $\mu\ge h^{-1}|\log h|^{-J}$. Namely the source of term containing logarithmic factor is the only perturbation of the type $C\mu^{-2}x_1^{-1}B(x',\mu^{-1}hD)$ but then if $\psi$ is even with respect to $x_1$ the results of calculation will be $0$ and if $\psi$ contains factor $x_1$ it would compensate $x_1^{-1}$ and no logarithm would appear.
\end{proof}

\subsubsection{}\label{sect-4-2-7} Now I need to consider the vicinity of $\Sigma$. Let us start from the general case first. 

Then scaling $x\to x/\gamma$, $h\to h/\gamma$, $\mu \to \mu\gamma$ and applying 
the same arguments as in the resonance case I estimate the contribution of zone
$\{\gamma (x)\asymp \gamma\}$ to the correction by 
\begin{equation*}
C\mu h^{-3} \gamma^{-1}\times (\mu^{-2}\gamma^{-1})^{1-\kappa}\times \gamma^2
\end{equation*}
which in comparison to resonance case gains an extra factor $\gamma$ and thus sums to its values as $\gamma=1$, which is $C\mu^{-1+\kappa} h^{-3}$ which in turn is $O(\mu^2h^{-2})$ as $\mu \ge h^{-1/3-\delta'}$. I ignore again term which sums to $C\mu^2h^{-2}$ in the end of the day.

Meanwhile contribution of zone $\{\gamma(x)\le {\bar\gamma}\}$ does not exceed 
$C\mu h^{-3}{\bar\gamma}^2= O(\mu ^2h^{-2})$. 

Thus I arrive to   the following statement:

\begin{proposition}\label{prop-4-21} Statement of proposition \ref{prop-4-3} remain true in the vicinity of $\Sigma$.
\end{proposition}

\begin{remark}\label{rem-4-22}  (i) Under condition 
\begin{equation}
\sum_{|\beta|\le K}|\nabla_\Sigma^\beta Vf_1^{-1}|\ge \epsilon_0\qquad \forall x
\label{4-32}
\end{equation}
with arbitrarily large $K$ one can skip a correction term  without deteriorating remainder estimate unless $h^{-1/3}|\log h|^{-J}\le \mu \le h^{-1/3}|\log h|^J$ and with the remainder estimate $O(h^{-8/3}|\log h|^J)$ in this border case.
I suspect that one can get rid off logarithmic factors and to prove an  estimate $O(\mu^{-1}h^{-3}+\mu^2h^{-2})$ even in the border case.

\medskip
\noindent
(ii) One can easily construct a correction term in the case $\Sigma\ne\emptyset$ resonances but an expression seems  rather too complicated. So I leave it to the curious reader.
\end{remark}

\subsubsection{}\label{sect-4-2-8} Finally let us consider the generic case near $\Sigma$. Using the same arguments as before I arrive to

\begin{proposition}\label{prop-4-23}  Under non-degeneracy conditions $(\ref{3-84})$  and $(\ref{3-85})$ for $\psi$ supported in the vicinity of $\Sigma$, the remainder estimate is $O(\mu^{-1}h^{-3})$ while the main part is given by the magnetic Weyl formula.
\end{proposition}

\section{Vanishing $V$ Case}
\label{sect-5}

\subsection{Generic Case}
\label{sect-5-1}
In the generic case, as $V=0$, $|\nabla V|\ge \epsilon$ and then the microhyperbolicity condition of \cite{IRO3} holds and then 

\begin{proposition}\label{prop-5-1}
 Assume that 
 \begin{equation}
 |V|+|\nabla V|\ge \epsilon.
 \label{5-1}
 \end{equation}
Then in  any dimension  without no condition to $F_{jk}$ other than 
$|\det (F_{jk})|\ge \epsilon$ for $\psi$ supported in the small vicinity of 
 \begin{equation}
\Delta \Def \{x:\ V(x)=0\}
 \label{5-2}
 \end{equation}
 the asymptotics  with the main part given by Magnetic Weyl expression and the remainder estimate $O(\mu^{-1}h^{1-d})$ holds.
\end{proposition}
This covers the generic case completely.

\subsection{General Case}
\label{sect-5-2}
The general case however is more complicated. Let us introduce a scaling function
\begin{equation}
\ell = \epsilon \bigl(|V|+|\nabla V|^2\bigr)^{1/2}+ {\frac 1 2}{\bar\ell}, \qquad
{\bar\ell} = \epsilon_0\max\bigl((\mu h)^{1/2},\mu^{-1}\bigr)
\label{5-3}
\end{equation}
with  dominating the second and the first terms in  ${\bar\ell}$ as $\mu\le h^{-1/3}$
and $\mu \ge h^{-1/3}$ respectively.

Let us apply scaling $x\mapsto x \ell$, $h\mapsto h\ell ^{-2}$, $\mu \mapsto \mu$, $V\mapsto V\ell^{-2}$: 

\begin{itemize}[label = - ]
\item Let  $\ell\le \epsilon_0(\mu h)^{1/2}$. Then I am in the classically forbidden zone and the contributions of $\ell$-element to the principal part and the remainder estimate are 0 and negligible respectively. 

\item Let $\ell\ge{\bar\ell}$ and $|\nabla V|\asymp \ell$. Then after rescaling I am in frames of subsection \ref{sect-5-1} and the contribution of $\ell$-element to the remainder is $O(\mu^{-1}h^{-3}\ell^6)$ and the total contribution of all 
$\ell$-elements to the remainder is $O(\mu^{-1}h^{-3}\ell^2)$ and I am done here. 

\item Also, the contribution of all elements with $\ell\asymp \mu^{-1}$ to the remainder is $O(\mu h^{-3}{\bar\ell}^2)$ which is $O(\mu^{-1}h^{-3})$. So case $\ell\asymp{\bar\ell}$ is completely covered.

\item Let us consider $V\asymp \ell^2$, $\ell\ge {\bar\ell}$; then condition $V\asymp 1$ is recovered after rescaling but now conditions to $F_{jk}$ could fail; actually these condition do not fail completely, but are replaced by somewhat weaker condition with extra factors $\ell$ or $\ell^2$  in the estimates from below. However in the weak magnetic field case this weakened condition is enough; I leave the details to the reader:
\end{itemize}

\begin{proposition}\label{prop-5-2}
As $F_{jk}$ is generic and $\mu \le h^{-\delta}$ the remainder estimate $O(\mu^{-1}h^{-3})$ holds.
\end{proposition}

Now one can assume that $h^{-\delta}\le \mu \le ch^{-1}$. Then since in section \ref{sect-3} no condition ``$|V|\ge \epsilon$'' was required, the remainder estimate is $O(\mu^{-1}h^{-3}+\mu^2h^{-2})$ but the principal part is given by the implicit formula (\ref{0-9}) rather than by the magnetic Weyl expression (\ref{4-4}) and now I need to modify arguments of section \ref{sect-4} to pass from (\ref{0-9}) to (\ref{4-4}).

Note that  in this implicit formula one always can take $T_0=\epsilon \mu^{-1}$ with an arbitrarily small constant $\epsilon >0$. Then in the virtue of arguments of the proof of propositions \ref{prop-4-3} , \ref{prop-4-4} and \ref{4-8}  the remainder estimate $O(\mu^{-1}h^{-3}+\mu ^2h^{-2-\delta})$ with an arbitrarily small exponent $\delta>0$. However in contrast to the analysis in subsection \ref{sect-4-1}, without condition ``$|V|\ge \epsilon$'' it does not translates into $T_0=Ch|\log h|$ but rather into $T_0= Ch\ell^{-2}|\log h|$ after the above partition is applied. However this last implicit formula (\ref{0-9})  translates (with the same error) into Weyl or Magnetic Weyl formula. Therefore  \emph{ theorem \ref{thm-0-2} with correction term $0$ is proven for
$\mu \le h^{-1/3-\delta}$ with an arbitrarily small exponent $\delta>0$\/}.

The similar arguments work as $\mu \ge h^{-1/3-\delta'}$ since one again refers to any $T_0\le \epsilon\mu^{-1}$.
So,

\begin{theorem}\label{thm-5-3} Both theorem \ref{thm-0-2} and  theorem \ref{thm-0-3}
\end{theorem}

\bibliographystyle{alpha}

\providecommand{\bysame}{\leavevmode\hbox to3em{\hrulefill}\thinspace}

\vglue .06truein

\begin{tabular}{rrl}
&{\hskip 220 pt} &Department of Mathematics,\cr
&&University of Toronto,\cr
&&40, St.George Str.,\cr
&&Toronto, Ontario M5S 2E4\cr
&&Canada\cr
&&ivrii@math.toronto.edu\cr
&&Fax: (416)978-4107\cr
\end{tabular}

\end{document}